\documentclass{amsart}
\usepackage[british,ngerman]{babel}
\usepackage{amsmath,amssymb}

\newcommand{\N}{\mathbb N}
\newcommand{\Z}{\mathbb Z}

\newcommand{\lcm}{\operatorname{lcm}}
\newcommand{\ggT}{\operatorname{ggT}}
\newcommand{\kgV}{\operatorname{kgV}}

\newcommand{\lra}{\longrightarrow}

\def\proofname{Beweis.}
\newtheorem{thm}{Satz}
\newtheorem{prop}[thm]{Proposition}
\newtheorem{lem}[thm]{Lemma}

\newenvironment{bew}{\begin{proof}}{\end{proof}}

\renewcommand{\theHsection}{de.\arabic{section}}

\begin{document}

\title{Zur Zahlentheorie der Griechen}
\author{Franz Lemmermeyer}

\begin{abstract}
In diesem Artikel stellen wir Euklids Fundamentalsatz der
Zahlentheorie vor und bringen Belege f\"ur die Behauptung,
dass dieser heute weitgehend unbekannte Satz den
Zahlentheoretikern vor Gau\ss{} so vertraut war wie uns die
Gau\ss{}sche Version des Fundamentalsatzes.

In einem zweiten Teil werden wir die Bedeutung von
Euklids Fundamentalsatz f\"ur die heutige Algebra erl\"autern.
Dabei gehen wir von Ringen aus, in denen das Gau\ss{}sche
Lemma gilt, und sto\ss{}en so auf Begriffe wie den der
ganzen Abgeschlossenheit und auf Dedekinds Prager Satz.
\end{abstract}

\maketitle

\begin{center}
{\Large Teil I. Euklids Fundamentalsatz der Arithmetik}
\end{center}
\medskip

%\section*{Teil I. Euklids Fundamentalsatz der Arithmetik}

Der Fundamentalsatz der Arithmetik, die Eindeutigkeit der
Primfaktorzerlegung in den nat\"urlichen Zahlen $\N$, wurde
erstmals von Gau{\ss} in seinen Disquisitiones Arithmeticae \cite{DA}
explizit formuliert und bewiesen. Bisweilen liest man, dieser
Satz gehe auf Euklid zur\"uck. Ganz falsch ist das sicher
nicht, besser allerdings w\"are es, die Frage etwas anders
zu stellen: es geht ja in erster Linie nicht darum, ob Euklid
wusste, dass sich jede Zahl eindeutig in ihre Primfaktoren zerlegen
l\"asst: den meisten Autoren, die sich z.B. mit der Konstruktion
vollkommener und befreundeter Zahlen befassten, waren zumindest
Spezialf\"alle dieses Resultats gel\"aufig, und
manche (nach Euklid insbesondere al-Farisi, sh.~\cite{AgFl},
und Prestet, sh.~\cite{Gold}) gaben daf\"ur auch Beweise \`a la
Euklid. Aber erst Gau\ss{} machte die eindeutige Primzerlegung zu
einem Fundamentalsatz.

Interessant ist daher vielmehr, ob Euklid dieses Ergebnis so
verstanden hat wie wir, d.h. ob er dessen fundamentalen Charakter
erkannt hat. Letzteres war nun zweifellos nicht der Fall. Dies
wirft dann die Frage auf, was denn die Grundlage der Euklidischen
Zahlentheorie war, und die Antwort darauf wird meist durch den Hinweis
auf den fundamentalen Charakter des euklidischen Algorithmus erledigt.
Die Hauptergebnisse der euklidischen Zahlentheorie, insbesondere das
Euklidische Lemma, wonach irreduzible Elemente prim sind, werden bei
Euklid aber alle mit Hilfe der Proposition VII.19 bewiesen. Wir wollen
im Folgenden zeigen, dass diese nur wenig bekannte Proposition nicht
nur f\"ur Euklid fundamental war, sondern von Euler und anderen
Zahlentheoretikern als ``Ersatz'' f\"ur die eindeutige Primzerlegung
benutzt wurde und selbst heute noch, wenn auch unter anderem Namen,
eine gewisse Rolle in der modernen Algebra spielt.

\section{Eindeutige Primzerlegung?}

Wenn man liest, dass Euklid die eindeutige Primzerlegung gekannt
und bewiesen habe, bezieht sich diese Aussage auf Euklids (vgl. die
Euklidausgaben von Thaer \cite{Thaer} und Heath \cite{Heath}):
\begin{quote}
  \emph{Proposition IX.14. Ist eine Zahl die kleinste, welche von
    gegebenen Primzahlen geteilt wird, dann wird sie von keinen
    anderen Primzahlen geteilt au{\ss}er den urspr\"unglichen.}
\end{quote}
Mit anderen Worten: ist $n = \kgV(p_1, \ldots, p_t)$ das kleinste
gemeinsame Vielfache\footnote{Dass das kgV teilerfremder Zahlen
gleich ihrem Produkt ist, folgt aus Euklids Proposition VII.34.}
der $p_j$, so sind $p_1$, \ldots, $p_t$ die einzigen Primfaktoren
von $n$. Anscheinend war es Heath \cite[vol. II, p. 402]{Heath}, der
darin die eindeutige Primzerlegung gesehen hat: in seinen
Kommentaren zu IX.14 schreibt er:
\begin{quote}
  \emph{In other words, a number can be resolved into prime factors
    in only one way.\footnote{Mit anderen Worten: eine Zahl kann
      auf nur eine Art und Weise in Primfaktoren zerlegt werden.}}
\end{quote}

Dass Euklid nicht das Produkt, sondern die kleinste durch
$p_1$, \ldots, $p_t$ teilbare Zahl benutzt, liegt daran,
dass Euklids geometrische Sprache ihm nur erlaubt hat, von
Produkten von zwei (ebene Zahlen) oder drei Zahlen (solide Zahlen)
zu sprechen.

Dass er den Beweis exemplarisch an einer Zahl, die durch
drei verschiedene Primzahlen teilbar ist, durchf\"uhrt, entspricht
den Gepflogenheiten in der griechischen Mathematik und ist eine
Konsequenz der Tatsache, dass den Griechen die Technik der
vollst\"andigen Induktion fehlte.

Die wesentliche Schw\"ache von IX.14 im Vergleich zum Satz
von der eindeutigen Primzerlegung ist die Tatsache, dass IX.14
die eindeutige Primzerlegung nur f\"ur Produkte verschiedener
Primzahlen liefert, also f\"ur quadratfreie Zahlen; letztendlich
besagt IX.14 von einem Produkt $n = p_1^{a_1} p_2^{a_2} p_3^{a_3}$
nur, dass $p_1$, $p_2$ und $p_3$ die einzigen Primfaktoren von $n$
sind\footnote{Es sei aber hinzugef\"ugt, dass Euklid an anderen
Stellen S\"atze beweist, mit denen sich IX.14 pr\"azisieren lie{\ss}e:
so wird in Buch IX u.a. bewiesen (IX.13), dass die einzigen Teiler einer
Primzahlpotenz $p^n$ die Zahlen $1$, $p$, $p^2$, \ldots, $p^{n-1}$
sind (f\"ur Euklid teilt keine Zahl sich selbst). Schlie\ss{}lich
hat Euklid auch bewiesen (IX.36), dass die einzigen Teiler der Zahl
$2^{p-1}q$ f\"ur prime $q = 2^p-1$ durch $1$, $2$, $4$, \ldots,
$2^{p-1}$, $q$, $2q$, $4q$, \ldots, $2^{p-2}q$ gegeben sind.}.

Trotz dieses Mangels wird Euklids IX.14 ebenso wie heute der
Satz von der Eindeutigkeit der Primzerlegung mit ``Euklids
Lemma'' (Prop. VII.30) bewiesen, wonach aus $p \mid ab$ f\"ur
prime\footnote{Der Begriff der Primzahl bei Euklid stimmt mit dem
\"uberein, was man in der heutigen Algebra irreduzibel nennt. Die
Eigenschaft aus Euklids Lemma benutzt man dann zur Definition von
Primzahlen.} $p$ immer $p \mid a$ oder $p \mid b$ folgt. Was
Euklid also von einer vollst\"andigen Formulierung nebst Beweis
des Fundamentalsatzes der Arithmetik abgehalten hat, war zum
einen sicherlich die mangelhafte Sprache, derer er sich bedienen
musste, zum andern aber auch die Tatsache, dass IX.14 f\"ur Euklid
alles andere als fundamental war und nicht viel mehr als ein kleines
Lemma, das er beil\"aufig erw\"ahnt. Euklid benutzt seine
Proposition IX.14 in der Folge an keiner Stelle, und diese Beobachtung
allein verbietet es, in diesem Zusammenhang von einem Fundamentalsatz
zu sprechen.

Das Fundament der euklidischen Zahlentheorie ist die Existenz des
euklidischen Algorithmus, und der Fundamentalsatz, aus dem die
zentralen S\"atze gewonnen werden, ist Euklids Proposition VII.19,
welcher wir uns nun zuwenden wollen.

\section{Euklids Proposition VII.19}

Um Euklids Proposition VII.19 verstehen zu k{\"o}nnen, ben{\"o}tigen
wir seine Definition VII.20:
\begin{quote}
  \emph{Zahlen sind proportional, wenn die erste Zahl dasselbe Vielfache,
    oder derselbe Teil, oder dieselben Teile, der zweiten ist wie
    die dritte der vierten.}
\end{quote}
Dass die erste Zahl dasselbe Vielfache der zweiten ist wie die
dritte der vierten, bedeutet, dass es eine nat\"urliche Zahl $m$
gibt mit $a = mb$ und $c = md$. Entsprechend bedeutet ``derselbe
Teil'', dass es einen Teiler $n$ von $b$ und $d$ gibt mit
$a = \frac1n b$ und $c = \frac1n d$; mit $b = nx$ und $d = ny$
liefert dies $a = x$, $b = nx$, $c = y$, $d = ny$. Endlich bedeutet
``dieselben Teile'' $a = \frac mn b$ und $c = \frac mn d$ f\"ur
nat\"urliche Zahlen $m, n$ mit $n \mid b$ und $n \mid d$.\footnote{Die
ersten beiden F\"alle w\"urden wir als die Spezialf\"alle $n = 1$
bzw. $m = 1$ der letzten Definition ansehen; Euklids Zahlbegriff
unterscheidet sich aber von unserem: f\"ur ihn waren Zahlen echte
Vielfache der Einheit $1$, und die Einheit keine Zahl.}
\"Ubersetzt bedeutet dies also:

\medskip\noindent{\bf Definition.}
F\"ur nat\"urliche Zahlen $a, b, c, d \in \N$ sind die
Verh\"altnisse $a:b$ und $c:d$ proportional (d.h. es gilt
$a:b = c:d$) genau dann, wenn es $m, n, x, y \in \N$ gibt mit
\begin{equation}\label{Eamx}
 a = mx, \quad b = nx, \quad c = my, \quad d = ny.
\end{equation}

So ist z.B. $6:2 = 9:3$ via $m = 3$, $n = 1$, $x = 2$ und $y = 3$.
Damit lautet
\begin{quote}
  \emph{Proposition VII.19. Sind vier Zahlen proportional, dann ist
    das Produkt der ersten und der vierten gleich dem Produkt der
    zweiten und der dritten Zahl. Ist das Produkt der ersten und der
    vierten gleich dem Produkt der zweiten und der dritten Zahl, dann
    sind die vier Zahlen proportional.}
\end{quote}
Euklid behauptet also, dass genau dann $a:b = c:d$ ist, wenn $ad = bc$
gilt. In der darauf folgenden Proposition VII.20 zeigt Euklid,
dass die minimalen Zahlen $a, b$ mit Verh\"altnis $c:d$ diese Zahlen
$c$ und $d$ teilen, und zwar so, dass $c = am$ und $d = bm$ gilt.
In VII.21 und VII.22 schlie\ss{}lich wird gezeigt, dass $a$ und $b$
genau dann die minimalen Zahlen mit $a:b = c:d$ sind, wenn $a$ und
$b$ teilerfremd sind.

Wir neigen dazu, Euklids Proportionen von Zahlen mit den modernen
Br\"uchen zu identifizieren (f\"ur Euklid waren Proportionen keine
Zahlen, und er erkl\"art auch nicht, wie man Proportionen addiert);
in dieser Sprache geht es in den Propositionen VII.19\,ff.\ darum zu
zeigen, dass sich jeder Bruch eindeutig in maximal gek\"urzter Form
darstellen l\"asst.

\subsection*{Von Zahlen und Gr{\"o}{\ss}en}
Neben Proportionen von Zahlen behandelt Euklid in Buch V auch
Proportionen von Gr{\"o}{\ss}en; darunter fallen Winkel, L\"angen,
Fl\"achen und Volumina, also Objekte, die sich kontinuierlich
\"andern k{\"o}nnen. Alle gleichartigen Gr{\"o}{\ss}en bei Euklid
haben folgende Eigenschaften: man kann sie addieren, und sie
verhalten sich archimedisch in folgendem Sinne: sind $A$ und $B$
Gr{\"o}{\ss}en derselben Art (solche lassen sich in Bezug auf ihre
Gr{\"o}{\ss}e vergleichen), so gibt es immer ein $n \in \N$ derart,
dass $nB > A$ ist.

F\"ur Verh\"altnisse $A:B$ und $C:D$ von Gr{\"o}{\ss}en\footnote{Hier
m\"ussen nat\"urlich $A$ und $B$ bzw. $C$ und $D$ gleichartige
Gr\"o\ss{}en sein.} gilt nach einer Eudoxos zugeschriebenen
Definition (Def.~V.5) $A:B = C:D$ genau dann, wenn gilt: f\"ur alle
nat\"urlichen Zahlen $m, n$ ist $mA < nB$ ($mA = nB$; $mA > nB$) genau
dann, wenn $mC < nD$ ($mC = nD$; $mC > nD$) ist.

Damit l\"asst sich dann z.B. zeigen, dass sich Kreisen
einbeschriebene \"ahnliche Polygone zueinander verhalten wie
die Quadrate \"uber den Kreisdurchmessern (Prop. XII.1), oder dass
sich Kreise zueinander verhalten wie die Quadrate \"uber ihren
Durchmessern (Prop. XII.2).

Um die Unterschiede der Definition der Proportionalit\"at von Zahlen
und Gr{\"o}{\ss}en deutlich zu machen, stellen wir Folgendes fest:
\begin{itemize}
\item Euklids Definition der Proportionalit\"at von Zahlen ist nicht
      nur in $\N$, sondern in beliebigen multiplikativ geschriebenen
      Monoiden sinnvoll.
\item Die Definition der Proportionalit\"at von Gr{\"o}{\ss}en ist
      f\"ur alle additiv geschriebenen und im obigen Sinne
      ``archimedischen'' Monoide sinnvoll.
\end{itemize}

Prinzipiell lie{\ss}en sich Zahlen auch als Gr{\"o}{\ss}en
auffassen\footnote{Denn $\N$ ist ja auch ein additives
archimedisches Monoid.}; Euklid macht das aber nicht. Die
manchmal anzutreffende Behauptung,
Euklid habe Aussagen \"uber Proportionen von Gr{\"o}{\ss}en bewusst
oder unbewusst auf Zahlen angewandt, scheint mir weitgehend
unbegr\"undet zu sein und passt auch nicht recht zur Tatsache, dass
die Lehre der Proportionen von Zahlen deutlich \"alter\footnote{Van
der Waerden schreibt in \cite{vdW} praktisch das komplette Buch VII
den Pythagoreern zu; wie wir weiter unten erl\"autern werden, dr\"angt
sich die Definition der Proportionalit\"at von Zahlen fast auf,
wenn man sich wie die Pythagoreer mit Harmonielehre und ganzen Zahlen
besch\"aftigt. Dagegen wird die Definition der Proportionalit\"at
von Gr\"o\ss{}en allgemein Eudoxos zugeschrieben.} ist als
diejenige von Gr{\"o}{\ss}en.

\section{Folgerungen aus Euklids VII.19}
Euklids VII.19 ist ein auf den ersten Blick unscheinbares Resultat;
tats\"achlich handelt es sich dabei aber um ein ganz
au{\ss}erordentliches Ergebnis. Im Folgenden wollen wir kurz auf
einige bekannte Folgerungen aus VII.19 eingehen. Den Anfang macht

\begin{prop}
Aus VII.19 folgt Euklids Lemma: alle irreduziblen Elemente sind prim.
\end{prop}

\begin{bew}
Sei $p$ irreduzibel und $p \mid bc$, also $pd = bc$. Nach VII.19
gibt es $m, n, x, y$ mit $p = mx$, $b = nx$, $c = my$ und $d = ny$.
Da $p$ irreduzibel ist, ist $m$ oder $x$ eine Einheit, und dies
impliziert $p \mid b$ im ersten und $p \mid c$ im zweiten Fall.
\end{bew}

VII.19 gen\"ugt auch zum Beweis der folgenden Verallgemeinerung von
Euklids Lemma:

\begin{prop}\label{P1}
Euklids Proposition VII.19 f\"ur nat\"urliche Zahlen ist \"aquivalent
mit dem folgenden Resultat: ist $a \mid bc$ und $\ggT(a,b) = 1$, dann
gilt $a \mid c$.
\end{prop}

\begin{bew}
Es gelte VII.19. Ist dann $a \mid bc$, also $ad = bc$ f\"ur ein
$d \in \N$, so gibt es $m, n, x, y \in \N$ mit $a = mx$, $b = nx$,
$c = my$ und $d = ny$. Wegen $\ggT(a,b) = 1$ muss daher $x = 1$
sein, folglich ist $a = m$ ein Teiler von $c = my$.

Sei nun $ad=bc$, und setze $x = \ggT(a,b)$. Dann gibt es
$m, n \in \N$ mit $a = mx$ und $b = nx$. Daher ist $md = nc$ und
$\ggT(m,n) = 1$; dies impliziert $m \mid c$, also $c = my$ f\"ur
ein $y \in \N$. Damit ist aber $d = \frac{nc}m = ny$, d.h. es gilt
(\ref{Eamx}).
\end{bew}

Es ist heutzutage nat\"urlich kein Problem, die Eindeutigkeit
der Primfaktorzerlegung mit Hilfe von Prop. \ref{P1} (oder dem
Euklidschen Lemma) zu beweisen. Unter Zuhilfenahme einer Idee
von Zermelo \cite{Zer} kann man das aber direkt aus Euklids VII.19
herleiten:

\begin{prop}
Die Eindeutigkeit der Zerlegung in irreduzible Elemente folgt aus
Euklids VII.19.
\end{prop}

\begin{bew}
W\"are die Eindeutigkeit falsch, so g\"abe es eine kleinste
nat\"urliche Zahl $N$ mit zwei verschiedenen Faktorisierungen.
Sei $p$ ein irreduzibler Faktor der ersten Faktorisierung, und
schreibe $N = ph$; entsprechend sei $q$ ein irreduzibler Faktor
der zweiten Faktorisierung, und $N = qk$. Aus $ph = qk$ folgt mit
VII.19 die Existenz von $m, n, x, y \in \N$ mit $p = mx$, $q = nx$,
$k = my$ und $h = ny$. Da $p$ irreduzibel ist, muss $m = 1$ oder
$x = 1$ sein. Dann folgt aber $p \mid q$ (und damit $p = q$) im
ersten und $p \mid k$ im zweiten Fall. Im zweiten Fall ist
$k < N$, nach Wahl von $N$ besitzt $k$ also eine eindeutige
Zerlegung in irreduzible Elemente; wegen $p \mid k$ kommt $p$ in
dieser vor, also auch in derjenigen Zerlegung von $k$, die aus der
zweiten Faktorisierung von $N$ durch Weglassen von $q$ entsteht.
In beiden F\"allen kommt $p$ also auch in der zweiten
Faktorisierung vor, und K\"urzen von $p$ liefert eine nat\"urliche
Zahl $M = \frac Np < N$ mit zwei verschiedenen Faktorisierungen in
irreduzible Elemente: dieser Widerspruch zur Annahme der
Minimalit\"at von $N$ beweist den Satz.
\end{bew}

Kalm\'ar \cite{Kalm} (vgl. Sur\'anyi \cite{Sur1,Sur2}) nennt Euklids
VII.19 den Vierzahlensatz und schreibt ihn Euler zu. Die
Verallgemeinerung auf mehr Zahlen (ausgehend von einer Faktorisierung
$a_1 \cdots a_i = b_1 \cdots b_j$) hat Euler ebenfalls schon gekannt;
sie folgt durch Induktion \"uber die Anzahl der Faktoren aus dem
einfachen Vierzahlensatz. Kalm\'ar hat dann bemerkt, dass aus dieser
Verallgemeinerung die Eindeutigkeit der Primzerlegung folgt (sh. auch
Bell \cite{Bel1,Bel2}, Rosenthall \cite{Ros} und Erd\H{o}s \& Sur\'anyi
\cite{ES}), und Sur\'anyi \cite[S. 44\,ff.]{Sur2} gab f\"ur den
Vierzahlensatz einen geometrischen Beweis. Es scheint bisher
niemandem aufgefallen zu sein, dass der Vierzahlensatz euklidischen
Ursprungs ist\footnote{Wer sich n\"aher mit den Elementen befasst
hat, h\"atte das nat\"urlich auf Anhieb sagen k{\"o}nnen, aber
anscheinend gibt es recht wenige Mathematiker, welche sowohl mit den
Elementen als auch mit der modernen Algebra hinreichend vertraut
sind. Selbst Sur\'anyi, der sich in \cite{Sur2} ganz explizit mit
griechischer Mathematik befasst, hat den Zusammenhang zwischen
Euklids VII.19 und dem Vierzahlensatz \"ubersehen.}.

Es ist eine leichte \"Ubung, z.B. mit Hilfe des Satzes von der
eindeutigen Primzerlegung das folgende Resultat zu beweisen, das
in vielen elementaren Untersuchungen \"uber diophantische Gleichungen
benutzt wird:

\begin{prop}\label{Psq}
Sind $a, b \in \N$ teilerfremd und ist $ab$ ein Quadrat, so sind
auch $a$ und $b$ Quadrate.
\end{prop}

Auch diese Proposition ist eine einfache Folgerung aus dem
Euklidischen VII.19, wie man leicht nachrechnet (wir werden
auf den Beweis in Prop. \ref{Psqc} noch einmal zur\"uckkommen).
Gau{\ss} hat in den Disquisitiones \cite{DA} als erster darauf
hingewiesen, dass dieses Ergebnis (bzw. die Verallgemeinerung
auf $n$-te Potenzen) aus der eindeutigen Primzerlegung in $\Z$
folgt. Der Umkehrschluss, dass sich Prop. \ref{Psq} nur mit
Hilfe des Fundamentalsatzes der Arithmetik beweisen l\"asst,
ist oft anzutreffen, aber falsch.

Euler, dem ja oft nachgesagt wird, er h\"atte die eindeutige
Primzerlegung in seinen Arbeiten unbewusst benutzt oder als
selbstverst\"andlich angesehen, hat Prop. \ref{Psq} in seiner
ersten Arbeit \cite{EuTQ} \"uber diophantische Fragen erw\"ahnt:
\begin{quote}
  \emph{Factum ex duobus pluribusque numeris inter se primis nec
    quadratum nec cubus nec ulla alia potestas esse potest, nisi
    singuli factores sint quadrata vel cubi vel eiusmodi aliae
    potestates.\footnote{Eine Zahl, die aus zwei oder mehr paarweise
      teilerfremden Faktoren besteht, kann weder ein Quadrat noch eine
      Kubikzahl noch irgendeine andere Potenz sein, au{\ss}er wenn die
      einzelnen Faktoren Quadrate, Kubikzahlen oder entsprechende
      andere Potenzen sind.}}
\end{quote}
Anstatt einen Beweis zu geben, bemerkt er dann:
\begin{quote}
  \emph{Demonstratio huius Lemmatis facilis est atque ab Euclide iam
    est tradita, ita ut superfluum foret eam hic
    exponere.\footnote{Der Beweis dieses Lemmas ist leicht und wurde
      bereits von Euklid \"uberliefert; es w\"are daher
      \"uberfl\"ussig, ihn hier darzulegen.}}
\end{quote}
Es ist mir nicht bekannt, auf welches Ergebnis von Euklid sich Euler
hier bezieht. F\"ur Quadrate folgt die Behauptung relativ einfach
aus Euklids Theorie der \"ahnlichen ebenen Zahlen (\"ahnliche ebene
Zahlen sind in heutiger Sprechweise Zahlen, deren Produkt ein Quadrat
ist; vgl. IX.1 und IX.2), und f\"ur h{\"o}here Potenzen mit etwas
mehr M\"uhe aus vielen seiner Propositionen im Buch VIII.

Wo wir uns also auf die eindeutige Primzerlegung, d.h. eine
einf\"uhrende Vorlesung in die Zahlentheorie oder die Gau\ss{}schen
Disquisitiones \cite{DA}, berufen w\"urden, erledigt Euler das mit
einem Verweis auf Euklid. Es scheint mir recht plausibel, dass Euler
auch an anderen Stellen, an denen wir eine L\"ucke sehen, die
durch eindeutige Primzerlegung geschlossen werden kann, etwaige
Einw\"ande mit einem Hinweis auf Euklid abgetan h\"atte.

In \cite{Eu152} zum Beispiel f\"uhrt er zum Studium befreundeter
Zahlen die Summe $\sigma(n)$ aller Teiler von $n \in \N$ ein
und zeigt in Lemma 1, dass $\sigma$ multiplikativ ist, also
$\sigma(mn) = \sigma(m)\sigma(n)$ f\"ur teilerfremde Zahlen $m, n$ gilt.
Hier ist sein Beweis: ``Denn das Produkt $mn$ hat erstens die
einzelnen Teiler beider Faktoren $m$ und $n$, und ist dann noch
teilbar durch die Produkte der Teiler von $m$ und derjenigen
von $n$. Deshalb ist die Summe aller Teiler von $mn$ gleich dem
Produkt von $\sigma(m)$ und $\sigma(n)$.'' Euler benutzt also das
folgende Lemma ohne Beweis:

\begin{lem}\label{LEu4}
Seien $m$ und $n$ teilerfremd. Dann l\"asst sich jeder Teiler $a$
von $mn$ eindeutig in der Form $a = rs$ schreiben, wo $r \mid m$
und $s \mid n$ gilt.
\end{lem}

Wir w\"urden dieses Lemma heute ohne zu z\"ogern aus dem
Fundamentalsatz der Arithmetik herleiten; Euler h\"atte dagegen
vermutlich auf Euklid verwiesen: Die Existenz von $r$ und $s$
folgt so: $a \mid mn$ bedeutet $ab = mn$ f\"ur ein $b \in \N$;
Euklids VII.19 liefert die Existenz von $r, s, t, u \in \N$ mit
$a = rs$, $b = tu$, $m = rt$ und $n = su$, und die Behauptung
folgt. Die Eindeutigkeit ist aber klar: ist $a = rs = r's'$ mit
$r, r' \mid m$ und $s, s' \mid n$, so ist $r \mid r's'$ und
$\ggT(r,s') \mid \ggT(m,n) = 1$, nach Prop. \ref{P1} also
$r \mid r'$. Ebenso folgt $r' \mid r$, also $r = r'$ und $s = s'$.

\section{Die L\"ucke}
Die G\"ultigkeit von Euklids VII.19 ist, wie wir gesehen haben,
f\"ur die Eindeutigkeit der Zerlegung in irreduzible Elemente
in $\N$ verantwortlich. Eine einfache Struktur, in welcher VII.19
nicht mehr richtig ist, ist das Monoid $M = \{1, 5, 9, 13, \ldots\}$
aller nat\"urlichen Zahlen $\equiv 1 \bmod 4$ (Hilbert hat in seiner
Vorlesung \cite{HZ} das Monoid aller nat\"urlichen Zahlen
$\equiv 1 \bmod 5$ benutzt, um die Einf\"uhrung von Idealen zu
motivieren). In $M$ gilt offenbar
$21 \cdot 21 = 9 \cdot 49$ (diese Gleichung zeigt auch sofort, dass
das Analogon von Prop. \ref{Psq} in $M$ nicht gilt), aber
$49:21 \ne 21:9$; denn aus (\ref{Eamx}) folgt die Existenz von
$m, n, x, y \in M$ mit $49 = mx$, $21 = nx = my$ und $9 = ny$. Da $9$
in $M$ irreduzibel ist, muss $n = 1$ oder $y = 1$ sein. Im ersten
Fall ist $y = 9$, also $9 \mid my = 21$, im zweiten Fall $n = 9$ und
damit ebenfalls $9 \mid nx = 21$; beides ist unm\"oglich.

Dagegen zeigt man leicht, dass f\"ur $a, b, f \in M$ immer
$a:b = af:bf$ gilt. Also ist
$$ 49 : 21 = 49\cdot 9 : 21 \cdot 9
           = 21 \cdot 21 : 21 \cdot 9 = 21 : 9, $$
wobei wir zweimal die Beobachtung $a:b = af:bf$ angewendet haben.
Wie wir aber schon gesehen haben, ist die Gleichung $49 : 21 = 21 : 9$
falsch. An welcher Stelle haben wir den Fehler gemacht? Die Gleichungen
$49:21 = 21 \cdot 21 : 21 \cdot 9$ und $21 \cdot 21 : 21 \cdot 9 = 21 : 9$
sind, wie man leicht auch direkt nachpr\"uft, zweifellos richtig. Der
Fehler in der obigen Argumentation liegt in der Annahme, dass die
Gleichheit von Proportionen\footnote{Man k\"onnte versucht sein, das
darauf zur\"uckzuf\"uhren, dass die Transitivit\"at der Gleichheit bei
Euklid ein Axiom ist (``Was demselben gleich ist, ist einander gleich'');
allerdings ist das, was wir Gleichheit von Proportionen nennen, keine
Gleichheit im Euklidischen Sinne: Euklid definiert nicht die
Gleichheit zweier Proportionen $a:b$ und $c:d$, sondern sagt,
vier Zahlen seien proportional, wenn es Zahlen wie in (\ref{Eamx})
gibt. Aus diesem Grund schreiben viele Kommentatoren Euklids nicht
$a:b = c:d$, sondern genauer $a:b :: c:d$.} transitiv ist: aus
$a:b = c:d$ und $c:d = e:f$ folgt im Allgemeinen nicht, dass auch
$a:b = e:f$ ist.

Taisbak hat in \cite{Tai} herausgefunden, dass Euklid im Beweis seines
Fundamentalsatzes VII.19 genau diese L\"ucke gelassen hat: er hat die
Transitivit\"at der Gleichheit benutzt, aber nicht bewiesen.
Wie diese L\"ucke zu schlie{\ss}en ist, kann man bei Taisbak
\cite{Tai} wie auch bei Pengelley \& Richman \cite{PR} nachlesen.
Das einfache Argument beruht darauf zu zeigen, dass man einige
der Zahlen im Vierzahlensatz mehr oder weniger ``kanonisch''
als gr{\"o}{\ss}te gemeinsame Teiler w\"ahlen kann:

\begin{lem}\label{Ltr}
Sei $a:b = c:d$. Dann gibt es $p, q, r, s \in \N$ mit $a = pr$,
$b = qr$, $c = ps$ und $d = qs$, wobei wir $r = \ggT(a,b)$ und
$s = \ggT(c,d)$ w\"ahlen k{\"o}nnen.
\end{lem}

\begin{bew}
Aus $a:b = c:d$ folgt die Existenz von Zahlen $m, n, x, y \in \N$
mit (\ref{Eamx}). Also ist $x \mid \ggT(a,b)$ und $y \mid \ggT(c,d)$,
folglich gibt es $i, j \in \N$ mit $\ggT(a,b) = ix$ und $\ggT(c,d) = jy$.
Aus $ix \mid a = mx$ und $ix \mid b = nx$ erhalten wir $i \mid m$ und
$i \mid n$, also $m = pi$ und $n = qi$ mit $p, q \in \N$. Insbesondere
ist $a = pix = p \cdot \ggT(a,b)$ und $b = q \cdot \ggT(a,b)$.

Jetzt behaupten wir $i \mid j$. Aus Symmetriegr\"unden (die Situation
ist symmetrisch unter Vertauschung von $(a,b,x)$ mit $(c,d,y)$) ist
dann auch $j \mid i$, also $i = j$. Aber dann bekommen wir wie
gew\"unscht $c = my = piy = pjy = p \cdot \ggT(c,d)$ und
$d = ny = qiy = qjy = q \cdot \ggT(c,d)$.

Es bleibt $i \mid j$ zu zeigen. Aus $m = pi$ und $n = qi$ folgt
$iy \mid my = c$ und $iy \mid ny = d$.
Also ist $iy \mid \ggT(c,d) = jy$, und damit $i \mid j$ wie behauptet.
\end{bew}

Damit k{\"o}nnen wir die Transitivit\"at der Gleichheit so zeigen:
sei $a:b = c:d$ und $c:d = e:f$. Nach Lemma \ref{Ltr} gibt es
$p, q \in \N$ mit $a = p \cdot \ggT(a,b)$, $b = q \cdot \ggT(a,b)$,
$c = p \cdot \ggT(c,d)$ und $d = q \cdot \ggT(c,d)$. Wenden wir
dieses Lemma ein zweites Mal an, so finden wir $p', q' \in \N$ mit
$c = p' \cdot \ggT(c,d)$, $d = q' \cdot \ggT(c,d)$,
$e = p' \cdot \ggT(e,f)$ und $f = q' \cdot \ggT(e,f)$. K\"urzen
liefert $p = p'$ und $q = q'$, und damit ist sicherlich $a:b = e:f$.

Man beachte, dass Lemma \ref{Ltr} in dem oben angef\"uhrten
Monoid $M$ nicht gilt, da in $M$ zwei Zahlen keinen ggT zu haben brauchen.
So sind die gemeinsamen Teiler von $3^2 \cdot 7 \cdot 11$ und
$3^2 \cdot 7 \cdot 19$ in $M$ gegeben durch $1$, $3^2$ und $3 \cdot 7$,
aber es gibt keinen gr{\"o}{\ss}ten gemeinsamen Teiler bez\"uglich der
Teilbarkeit, also keinen gemeinsamen Teiler, der durch alle gemeinsamen
Teiler teilbar ist.

\section{Euklid oder Pythagoras?}

Zahlreiche Mathematikhistoriker schreiben Euklids Zugang zur
Zahlentheorie den Pythagoreern zu. Diese sind, wie man annimmt,
auf dem Weg \"uber die Musik zu ihren Proportionen gekommen: sie
hatten beobachtet, dass die L\"angen schwingender Saiten in
einfachen Verh\"altnissen zueinander stehen, wenn die entsprechenden
T{\"o}ne harmonisch klingen. So entsprechen die Verh\"altnisse $1:2$
bzw. $2:3$ der Oktave bzw. der Quinte.

Eine der ersten Beobachtungen d\"urfte dann gewesen sein, dass sich
gewisse Proportionen vereinfachen lassen: so ist ja $2:4 = 1:2$.
Letztlich f\"uhrt diese Fragestellung zwangsl\"aufig auf das
Ergebnis, dass sich Proportionen eindeutig k\"urzen lassen: sind
$m$ und $n$ minimal mit $a:b = m:n$, so gilt $m \mid a$ und $n \mid b$;
genauer: es gibt ein $x \in \N$ mit $a = mx$ und $b = nx$.
Das ist aber der Inhalt von Euklids Prop. VII.20. Danach wird in
VII.22 gezeigt, dass die minimalen $m$ und $n$ teilerfremd sind.

Hat man VII.20 einmal empirisch entdeckt, und ist $a:b = c:d$, so folgt
einerseits $a = mx$ und $b = nx$, andererseits $c = my$ und $d = ny$.
Das ist aber gerade die Euklidische Definition der Proportionalit\"at
von Zahlen, die sich damit ganz nat\"urlich aus der Besch\"aftigung
mit elementarsten Eigenschaften von Proportionen ergibt.

Nach Meinung von van der Waerden \cite{vdW} machten die Pythagoreer
diesen Vierzahlensatz zur Grundlage ihrer Zahlentheorie, welche
Euklid dann sp\"ater samt Aufbau im Wesentlichen unver\"andert
\"ubernahm. Ebenso wie andere Kommentatoren vor und nach ihm hat
van der Waerden Euklids Definition der Proportionalit\"at von Zahlen
eher als Fehltritt\footnote{D. Fowler schreibt z.B. in \cite{Fow}:
\begin{quote}
  \emph{Commentators rarely point out the unsatisfactory nature of
    this definition: it is a vivid, though incomplete, description of
    four numbers in proportion, and not the mathematical criterion
    needed for the foundation of a theory.}
\end{quote}}
denn als Geniestreich betrachtet. Allerdings zeigen die obigen
Betrachtungen und z.B. der Beweis von Lemma \ref{LEu4}, dass die
Definition einerseits recht nat\"urlich und andererseits der
Vierzahlensatz zum Beweis fundamentaler Ergebnisse der diophantischen
Analysis durchaus geeignet ist.

Au{\ss}er den von van der Waerden angef\"uhrten Argumenten f\"ur die
Behauptung, Euklids Buch VII gehe auf die Pythagoreer zur\"uck,
kann man auch das Fehlen des Beweises der Transitivit\"at der
Gleichheit f\"ur Proportionen von Zahlen als Beleg daf\"ur auffassen:
bei Zugrundelegung der von Eudoxos gegebenen Definition der Gleichheit
von Proportionen wird der Beweis der Transitivit\"at der Gleichheit
fast trivial. Dennoch hat ihn Euklid in V.11 ausgef\"uhrt. Dies mag
daran liegen, dass die Lehre der Proportionen von Gr{\"o}{\ss}en
auf Eudoxos zur\"uckgeht, und dass zu dessen Zeit die Einsicht, dass
derartige Schritte bewiesen werden m\"ussen, unter den Mathematikern
bereits verbreitet war. Es ist meiner Ansicht nach durchaus denkbar,
dass die Pythagoreer schlicht \"ubersehen haben, dass die Transitivit\"at
der Gleichheit von Proportionen von Zahlen \"uberhaupt nicht trivial ist.
Dass auch Euklid das nicht gesehen hat, obwohl er nat\"urlich V.11 kannte,
verwundert schon eher, lie{\ss}e sich aber dadurch erkl\"aren, dass er die
Zahlentheorie der Pythagoreer ohne gro{\ss}e \"Anderungen in sein Buch
VII \"ubernommen hat.

\section{Der Vierzahlensatz}

Um uns \"uber die Bedeutung des Euklidischen Vierzahlensatzes
Klarheit zu verschaffen, nennen wir Integrit\"atsringe, in denen
dieser Satz gilt, Rieszsche\footnote{Zafrullah \cite{Z} hat
diese Ringe pr\"a-Schreiersch genannt; weiter unten werden wir
aber sehen, dass die Teilbarkeitsgruppen Rieszscher Ringe ``Rieszsche
Gruppen'' im Sinne der Gruppentheorie sind.} Ringe.
In einem Rieszschen Ring $R$ folgt also aus $ad = bc$ mit
$a, b, c, d \in R$ die Existenz von $m, n, x, y \in R$ mit
$a = mx$, $b = nx$, $c = my$ und $d = ny$.

Ein Element $a \in R$ hei{\ss}t primal\footnote{Dieser Begriff geht
auf Cohn \cite{Cohn} zur\"uck.}, wenn aus $a \mid bc$ f\"ur beliebige
$b, c \in R$ die Existenz von $r, s \in R$ folgt mit $a = rs$, $r \mid b$
und $s \mid c$. Man beachte, dass $0$ und die Einheiten immer primal sind.

Ein Integrit\"atsring ist offenbar genau dann Rieszsch, wenn jedes
Element von $R$ primal ist. Die folgende Proposition zeigt, dass in
Rieszschen Ringen alle irreduziblen Elemente prim sind:

\begin{prop}\label{LIP}
Irreduzible Elemente sind genau dann primal, wenn sie prim sind.
\end{prop}

\begin{bew}
Sei $a$ prim und $a \mid cd$; dann ist o.B.d.A. $a \mid c$. Mit $r = a$
und $s = 1$ gilt dann $a = rs$, $r \mid c$ und $s \mid d$, d.h.
$a$ ist primal.

Sei nun $a$ irreduzibel und primal. Ist $a \mid cd$, so gibt es
$r, s$ mit $a = rs$, $r \mid c$ und $s \mid d$. Da $a$ irreduzibel
ist, muss $r$ oder $s$ -- sagen wir $r$ -- eine Einheit sein.
Aus $s \mid d$ folgt dann aber $a \mid d$, und damit ist $a$ prim.
\end{bew}

Euklid hat in Prop. VII.19 im Wesentlichen gezeigt, dass $\Z$ Rieszsch
ist; sein Beweis, dass irreduzible Elemente in $\N$ prim sind, ist im
Wesentlichen eine H\"alfte des obigen Beweises.

Sieht man sich den Euklidischen Beweis von VII.19 (bzw. dessen
Erg\"anzung durch Taisbak) genauer an, so stellt man fest, dass
der Vierzahlensatz in allen GGT-Ringen (dies sind Integrit\"atsringe,
in denen je zwei Elemente einen gr{\"o}{\ss}ten gemeinsamen Teiler
besitzen) gilt, dass also jeder GGT-Ring Rieszsch ist.

Es sei an dieser Stelle darauf hingewiesen, dass primale Elemente
auf der einen und irreduzible bzw. prime Elemente auf der andern
Seite von ganz unterschiedlicher Natur sind: so sind Produkte von
mehreren irreduziblen (bzw. primen) Elementen nicht mehr irreduzibel
(bzw. prim), abgesehen vom trivialen Fall eines einzigen Faktors.
F\"ur primale Elemente dagegen gilt:

\begin{lem}
Produkte primaler Elemente sind primal.
\end{lem}

\begin{bew}
Seien $a, b \in R$ primal; ist $a = 0$ oder $b = 0$, so ist $ab = 0$
primal, wir d\"urfen also $a, b \ne 0$ annehmen. Sei $ab \mid c_1c_2$,
also $c_1c_2 = abd$.
Da $a$ primal ist, gibt es eine Faktorisierung $a = a_1a_2$ mit
$a_i \mid c_i$. Schreiben wir $c_i = a_id_i$, so folgt
$abd = c_1c_2 = a_1a_2d_1d_2 = ad_1d_2$, also $bd = d_1d_2$.
Da $b$ primal ist, gilt $b = b_1b_2$ mit $b_i \mid d_i$.

Damit ist $ab = a_1b_1a_2b_2$ mit $a_ib_i \mid a_id_i = c_i$,
und folglich ist $ab$ primal. Mit Induktion folgt dann, dass
beliebige Produkte primaler Elemente primal sind.
\end{bew}

Jetzt zeigen wir folgende Verallgemeinerung von Prop. \ref{Psq}:

\begin{prop}\label{Psqc}
Sei $R$ Rieszsch. Gilt dann $ab=c^2$ f\"ur $a, b, c \in R$ und sind
$a$ und $b$ teilerfremd, dann gibt es eine Einheit $e \in R^\times$
derart, dass $ae$ und $be$ Quadrate in $R$ sind.
\end{prop}

\medskip\noindent
{\bf Bemerkung.}
Obwohl der Begriff eines gr{\"o}{\ss}ten gemeinsamen Teilers
zweier beliebiger Ringelemente nur in GGT-Ringen sinnvoll ist,
darf man in beliebigen Integrit\"atsringen von ``teilerfremden''
Elementen reden: dies sind per definitionem solche, deren gemeinsame
Teiler allesamt Einheiten sind. So sind z.B. $2$ und $3$ teilerfremd
in $R = \Z[\sqrt{-5}\,]$, w\"ahrend $6$ und $2+2\sqrt{-5}$ in $R$
keinen ggT besitzen.

\medskip

\begin{proof}[Beweis von Prop. \ref{Psqc}]
Nach dem Vierzahlensatz gibt es $m, n, x, y \in R$ mit
$a = mx$, $b = ny$ und $c = my = nx$. Wendet man den
Vierzahlensatz auf die Gleichung $my = nx$ an, so folgt
die Existenz von $r, s, t, u \in R$ mit $m = rs$, $y = tu$,
$n = rt$ und $x = su$. Mit $a$ und $b$ sind nat\"urlich
auch $m$ und $n$ bzw. $x$ und $y$ teilerfremd, folglich sind
$r$ und $u$ Einheiten in $R$. Mit $e = 1/ru$ ist dann
$ae = mxe = rs^2ue = s^2$ und $be = nye = rt^2ue = t^2$.
\end{proof}

Auch wenn sich einige von Euler ohne Beweis verwendete Tatsachen
mit Euklids Vierzahlensatz begr\"unden lassen, sind Rieszsche
Ringe doch nur ``fast'' faktorielle Ringe. Betrachten wir zwei
Faktorisierungen $a_1 \cdots a_m = b_1 \cdots b_n$, so folgt in
Rieszschen Ringen nur, dass beide Faktorisierungen eine gemeinsame
Verfeinerung besitzen: es existieren Elemente $c_{ij} \in R$
($1 \le i \le m$, $1 \le j \le n$) mit $a_i = c_{i1} \cdots c_{in}$
und $b_j = c_{1j} \cdots c_{mj}$. F\"ur $m = n = 2$ ist dies gerade
der Vierzahlensatz, der allgemeine Fall (der durch Induktion aus dem
Vierzahlensatz folgt) hei\ss{}t die ``Schreiersche
Verfeinerungseigenschaft''. Beispielsweise hat der Ring aller ganzen
algebraischen Zahlen diese Verfeinerungseigenschaft, ist aber kein
faktorieller Ring mangels irreduzibler Elemente.

Integrit\"atsringe, in denen sich jedes von $0$ verschiedene Element
als Produkt einer Einheit und irreduzibler Elemente schreiben l\"asst,
nennt man atomar. Damit gilt der folgende

\begin{thm}
Ein Integrit\"atsring $R$ ist genau dann faktoriell, wenn $R$
Rieszsch und atomar ist.
\end{thm}

\begin{bew}
In Rieszschen Ringen sind irreduzible Elemente prim nach Prop.
\ref{LIP}; damit sind Faktorisierungen in irreduzible Elemente
bekanntlich eindeutig. Ist umgekehrt $R$ faktoriell, so ist $R$
ein GGT-Ring und damit Rieszsch, und zweitens hat jedes Element
$\ne 0$ eine Faktorisierung in irreduzible (sogar prime) Elemente.
\end{bew}

Diese Aussage ist (wegen Prop. \ref{LIP}) etwas schw\"acher als
der bekannte Satz, wonach Integrit\"atsringe genau dann faktoriell
sind, wenn sie atomar sind und jedes irreduzible Element prim ist.

Der Euklidische Aufbau der Zahlentheorie kann wie folgt
schematisiert werden:
$$ \text{Euklid} \lra \text{GGT} \lra \text{Riesz} \lra \text{AP}, $$
wobei AP f\"ur die G\"ultigkeit des Euklidischen Lemmas steht, also
f\"ur die Eigenschaft, dass irreduzible Elemente (Atome) prim sind.
Allgemeiner kann man folgende Inklusionen zeigen:
\begin{figure}[h!]
$$ \begin{array}{cccccccccccc}
\text{Euklid} & \subset &  \text{HI}    & \subset & \text{Bezout}
        & & & & & & \\
& &     \cap      &         &   \cap     &         & & & & & \\
& & \text{ZPE} & \subset & \text{GGT} & \subset & \text{Riesz}
              & \subset & \text{AP} \\
& &     \cap      & &       &    & & & & & \\
& & \text{atomar}  & &       &    & & & & &
\end{array}$$
%\caption{Klassen von Ringen}\label{FR}
\end{figure}

Hierbei bezeichnen HI Hauptidealringe, ZPE faktorielle Ringe (also
solche mit eindeutiger Zerlegung in Prim-Elemente), und AP Ringe, in
denen alle Atome (also alle irreduziblen Elemente) prim sind.
Bezoutringe sind solche, in denen je zwei Elemente $a, b$ einen ggT
besitzen und sich dieser als $R$-Linearkombination von $a$ und $b$
schreiben l\"asst; gleichbedeutend damit ist, dass die Summe
zweier Hauptideale wieder ein Hauptideal ist.

W\"ahrend Hilbert die euklidische Geometrie auf eine moderne
axiomatische Grundlage gestellt hat, ist eine entsprechende
Untersuchung der euklidischen Zahlentheorie meines Wissens
bisher nicht durchgef\"uhrt worden. Man kann sich beispielsweise
fragen, in welchen Ringen Folgendes gilt: sind $a$ und $b$
teilerfremd (also jeder gemeinsame Teiler eine Einheit) und ist
$a \mid bc$, dann gilt $a \mid c$. Der erste Teil des Beweises
von Prop. \ref{P1} zeigt, dass dies in Rieszschen Ringen gilt;
der zweite Teil des Beweises benutzt aber die Existenz des ggT,
sodass wir nicht schlie{\ss}en k{\"o}nnen, dass solche Ringe
auch Rieszsch sind. Andererseits sind Ringe mit dieser
Eigenschaft AP. Man kann nun allgemein die Frage stellen, in
welchen Integrit\"atsringen die Propositionen aus Euklids
Zahlentheorie gelten, d.h. welche Axiome man zu ihrem Beweis
ben\"otigt.

\section{Rieszsche Gruppen}

In diesem Abschnitt wollen wir auf einige Verallgemeinerungen des
Vierzahlensatzes eingehen. Dazu werden wir etwas allgemeiner
Rieszsche Gruppen untersuchen.

Sei $R$ ein Integrit\"atsring mit Quotientenk{\"o}rper $K$. Wir
erinnern daran, dass $a \mid c$ f\"ur $a, c \in K^\times$ bedeutet,
dass es ein $b \in R$ gibt mit $ab = c$. F\"ur die Untersuchung von
Teilbarkeitseigenschaften sind Einheiten irrelevant; es ist
daher sinnvoll, statt $K^\times$ die Faktorgruppe $G(R) = K^\times/R^\times$
zu betrachten, die wir auch als Gruppe der von $(0)$ verschiedenen
gebrochenen Hauptideale auffassen k{\"o}nnen; die Elemente von $G(R)$
bezeichnen wir vorl\"aufig mit $aR^\times = (a)$. Die Gruppe $G(R)$ ist
halbgeordnet durch $(a) \le (b)$ genau dann, wenn $a \mid b$, d.h.
wenn $\frac ba \in R$ ist. Weiter ist $(a) \le (b)$ \"aquivalent mit
$(ac) \le (bc)$ f\"ur ein beliebiges $(c) \in G(R)$. Die Menge
$G(R)^+ = \{(g) \in G(R): (g) \ge (1)\}$ aller ``positiven'' Elemente
besteht also aus allen ganzen Hauptidealen $\ne (0)$ und ist offenbar
ein Monoid, das aus $R^* = R \setminus \{0\}$ durch \"Ubergang zu
Idealen, also durch Identifizieren assoziierter Elemente entsteht.

Wir wollen nun den Euklidischen Vierzahlensatz f\"ur $R$ in diese
Sprache \"ubersetzen. Ist $ad = bc$, so gilt $a \mid bc$, folglich
$(a) \le (b)(c)$ f\"ur die entsprechenden Hauptideale in $G(R)$.
Nach dem Vierzahlensatz existieren $r, s, t, u \in R$ mit $a = rs$,
$b = ru$, $c = ts$ und $d = tu$. Also gilt $r \mid b$ und $s \mid c$,
und damit $(r) \le (b)$ sowie $(s) \le (c)$. Daher gilt

\begin{prop}
Sei $R$ ein Integrit\"atsring. Dann gilt der Vierzahlensatz in $R$ genau
dann, wenn $G(R)$ die folgende ``Fuchssche Interpolationseigenschaft''
$F_2$  besitzt\footnote{Ab jetzt bezeichnen wir die Elemente von $G(R)$
statt mit $(a)$ einfach mit $a$.}: Sind $a, b, c \in G(R)^+$ und ist
$a \le bc$, dann gibt es eine Faktorisierung $a = rs$ mit
$r, s \in G(R)^+$ derart, dass $r \le b$ und $s \le c$ ist.
\end{prop}

Wir werden nun eine Reihe allgemeiner Interpolations- und
Verfeinerungseigenschaften von halbgeordneten abelschen Gruppen $G$
definieren:
\begin{itemize}
\item Die ``Fuchssche Interpolationseigenschaft'' $F_n$: Sind
      $a$, $b_1$, \ldots, $b_n \in G^+$, und ist $a \le b_1 \cdots b_n$,
      so gibt es $a_i \in G^+$ mit $a = a_1 \cdots a_n$ und $a_i \le b_i$
      f\"ur $1 \le i \le n$.
\item Die ``Rieszsche Interpolationseigenschaft'' $R_{m,n}$:
      Sind $a_1, \ldots, a_m$, $b_1, \ldots, b_n \in G$ gegeben mit
      $a_i \le b_j$ f\"ur alle $1 \le i \le m$ und $1 \le j \le n$,
      dann gibt es ein $c \in G$ mit $a_i \le c \le b_j$ f\"ur alle
      $i$ und $j$.
\item Die ``Schreiersche Verfeinerungseigenschaft'' $S_{m,n}$: Sind
      $a_1, \ldots, a_m$, $b_1, \ldots, b_n \in G^+$ gegeben mit
      $a_1 \cdots a_m = b_1 \cdots b_n$, so gibt es Elemente
      $c_{ij} \in G^+$ ($1 \le i \le m$, $1 \le j \le n$) mit
      $a_i = c_{i1} \cdots c_{in}$ und $b_j = c_{1j} \cdots c_{mj}$.
\item Die ``Euklidische Verfeinerungseigenschaft'' $E_{m,n}$: Sind
      $a, b_i, c_j \in G^+$ und $a \le b_ic_j$ f\"ur $i = 1$, \ldots, $m$
      und $j = 1$, \ldots, $n$, dann gibt es eine Faktorisierung $a = rs$
      mit $r, s \in G^+$, $r \le b_i$ und $s \le c_j$ f\"ur alle
      $1 \le i \le m$ und $1 \le j \le n$.
\end{itemize}

Die Eigenschaften $R_{2,0}$ und $R_{0,2}$ kann man so auffassen: zu
$a_1, a_2 \in G$ gibt es ein $c \in G$ mit $a_1 \le c$ und $a_2 \le c$
(das ist $R_{2,0}$) bzw. ein $b \in G$ mit $b \le a_1$ und $b \le a_2$
(das ist $R_{0,2}$); halbgeordnete Gruppen mit dieser Eigenschaft nennt
man gerichtet\footnote{Da $a \le b$ zu $b^{-1} \le a^{-1}$ \"aquivalent
ist, ist eine halbgeordnete Gruppe bereits gerichtet, wenn z.B.
$R_{2,0}$ erf\"ullt ist.}. Die Teilbarkeitsgruppen $G(R)$ von
Integrit\"atsringen sind automatisch gerichtet: zu
$a_1, a_2 \in K^\times$ gibt es, wie man sofort sieht, immer ein
$c \in K^\times$ mit $a_1 \mid c$, $a_2 \mid c$ (man w\"ahle
$c$ als das Produkt der (beliebig gew\"ahlten) Z\"ahler von $a_1$
und $a_2$). Analog findet man immer ein $b \in K^\times$ mit
$b \mid a_1$ und $b \mid a_2$.

Jetzt behaupten wir

\begin{thm}\label{SVerf}
F\"ur eine halbgeordnete gerichtete abelsche Gruppe $G$ sind folgende
Eigenschaften gleichbedeutend:
\begin{enumerate}
\item\label{R1} In $G$ gilt $R_{2,2}$;
\item\label{R2} In $G$ gilt $R_{m,n}$ f\"ur alle $m, n \ge 0$;
\item\label{R3} In $G$ gilt $F_2$;
\item\label{R4} In $G$ gilt $F_n$ f\"ur alle $n \ge 1$;
\item\label{R5} In $G$ gilt $S_{2,2}$;
\item\label{R6} In $G$ gilt $S_{m,n}$ f\"ur alle $m, n \ge 1$;
\item\label{R7} In $G$ gilt $E_{1,1}$;
\item\label{R8} In $G$ gilt $E_{m,n}$  f\"ur alle $m, n \ge 1$;
\item\label{R9} In $G$ gilt $[1,a] \cdot [1,b] = [1,ab]$ f\"ur alle
       $a, b \in G^+$ (hier ist $[a,b] = \{r \in G: a \le r \le b\}$).
\end{enumerate}
\end{thm}

Die Voraussetzung, dass $G$ gerichtet ist, braucht man nur, um von
$R_{2,2}$ auf $R_{m,n}$ schlie\ss{}en zu k\"onnen. Die Beweise sind
alle elementar (im Wesentlichen vollst\"andige Induktion) und brauchen
hier nicht vorgef\"uhrt zu werden.

Halbgeordnete abelsche Gruppen, welche die in Satz \ref{SVerf}
angef\"uhrten Eigenschaften haben, nennt man nach L.~Fuchs \cite{Fu}
Rieszsche Gruppen\footnote{Genauer wird bei Rieszschen Gruppen auf
die Forderung der Kommutativit\"at verzichtet; in diesem Fall sind
nicht mehr alle Eigenschaften aus Satz \ref{SVerf} gleichbedeutend.}.
Diese wurden erstmals von F.~Riesz \cite{Rie} 1940 untersucht, danach
von Birkhoff \cite{Birkh} und H.~Bauer \cite{Bauer}. In diesen
Arbeiten findet man auch Teile des Beweises von Satz \ref{SVerf}
und weitere zu \ref{R1} -- \ref{R9} \"aquivalente Eigenschaften.
Die Schreiersche Verfeinerungseigenschaft $S_{m,n}$ ist die
Verallgemeinerung des Vierzahlensatzes, welche Kalm\'ar zum Beweis
des Fundamentalsatzes der Arithmetik benutzt hat. Im zweiten Teil
werden wir schlie\ss{}lich die Euklidische Verfeinerungseigenschaft
ben\"otigen.

\section*{Abschlie{\ss}ende Bemerkungen}

Wer sich n\"aher mit der Zahlentheorie Euklids befasst,
st\"o\ss{}t schnell auf die zentrale Rolle seiner Proposition
VII.19, des Vierzahlensatzes. Aber auch wer sich nur f\"ur die
Geschichte der Zahlentheorie vor Gau\ss{} interessiert, muss
Euklids zahlentheoretische B\"ucher, und insbesondere das zentrale
Ergebnis VII.19, kennen: der oft formulierte Vorwurf, Euler h\"atte
an manchen Stellen implizit ein Ergebnis von Gau\ss{}, n\"amlich die
eindeutige Primzerlegung, verwendet, ist in dieser Form erstens
ohnehin etwas anachronistisch und vor dem Hintergrund von Euklids
VII.19 zweitens wohl auch nicht haltbar\footnote{Auch andere
Autoren (sogar noch nach Gau\ss{}) haben den Vierzahlensatz benutzt,
und zwar ohne auf Euler oder Euklid zu verweisen: sh.
Barlow \cite[p. 401]{Bar11}, Kramer \cite{Kramer} oder auch
P\'epin \cite[p. 24]{PepLeg}, \cite[p. 413]{Pep752}.}.
Mit diesen Bemerkungen soll Euler nicht von jeglicher ``Schuld''
reingewaschen werden: eine erkleckliche Anzahl seiner Beweise sind
nach heutigem Verst\"andnis nicht rigoros, und nicht alle L\"ucken
lassen sich mit Euklid f\"ullen. Der Pauschalvorwurf aber, Euler
h\"atte oft stillschweigend die eindeutige Primzerlegung benutzt,
sollte vielleicht doch \"uberdacht werden.

\newpage

\begin{otherlanguage*}{british}

\renewcommand{\proofname}{Proof}
\renewcommand{\theHsection}{en.\arabic{section}}

\begin{center}
{\Large Part I. Euclid's Fundamental Theorem of Arithmetic}
\end{center}
\medskip

\medskip\noindent
\begin{quote}
  \begin{small}
\textsc{Abstract.}
In this article we present Euclid's fundamental theorem of number
theory and give evidence for the claim that this theorem, which is
largely unknown today, was as familiar to number theorists before
Gauss as the Gaussian version of the fundamental theorem is to us.

In a second part we will explain the significance of Euclid's
fundamental theorem for modern algebra. Starting from rings in
which Gauss's Lemma holds, we will be led to notions such as
integral closure and to Dedekind's Prague theorem.
  \end{small}  
\end{quote}

\bigskip

The fundamental theorem of arithmetic, i.e., the uniqueness of prime
factorization in the natural numbers $\N$, was first explicitly
formulated and proved by Gauss in his Disquisitiones Arithmeticae
\cite{DA}. Occasionally one reads that this theorem goes back to
Euclid. This is certainly not entirely wrong; it would be better,
however, to pose the question somewhat differently: the main point
is not whether Euclid knew that every number can be factored uniquely
into primes. Most authors who studied, for example, the construction
of perfect and amicable numbers were familiar with at least special
cases of this result, and some (after Euclid, in particular
al-Farisi, see \cite{AgFl}, and Prestet, see \cite{Gold}) also gave
proofs \`a la Euclid. But it was only Gauss who turned unique prime
factorization into a fundamental theorem.

What is interesting, therefore, is rather whether Euclid understood
this result the way we do, i.e., whether he recognized its
fundamental character. This was certainly not the case. This raises
the question of what the foundation of Euclidean number theory
actually was, and the answer is usually given by pointing to the
fundamental character of the Euclidean algorithm. In Euclid, however,
the main results of his number theory, in particular Euclid's Lemma,
according to which irreducible elements are prime, are all proved
with the help of Proposition VII.19. In what follows we want to show
that this little-known proposition was not only fundamental for
Euclid, but was used by Euler and other number theorists as a
``substitute'' for unique prime factorization, and that even today,
albeit under a different name, it still plays a certain role in
modern algebra.

\section{Unique prime factorization?}

When one reads that Euclid knew and proved unique prime
factorization, this statement refers to Euclid's (cf.\ the editions
of Euclid by Thaer \cite{Thaer} and Heath \cite{Heath}):
\begin{quote}
  \emph{Proposition IX.14. If a number be the least that is measured
    by prime numbers, it will not be measured by any other prime
    number except those originally measuring it.}
\end{quote}
In other words: if $n = \lcm(p_1, \ldots, p_t)$ is the least common
multiple\footnote{That the lcm of coprime numbers equals their
product follows from Euclid's Proposition VII.34.} of the $p_j$,
then $p_1$, \ldots, $p_t$ are the only prime factors of $n$.
Apparently it was Heath \cite[vol. II, p. 402]{Heath} who saw unique
prime factorization in this statement; in his commentary on IX.14 he
writes:
\begin{quote}
  \emph{In other words, a number can be resolved into prime factors in
    only one way.}
\end{quote}

The fact that Euclid uses not the product but the smallest number
divisible by $p_1$, \ldots, $p_t$ is due to the fact that Euclid's
geometric language only allowed him to speak of products of two
(plane numbers) or three numbers (solid numbers).

That he carries out the proof exemplarily for a number divisible by
three different primes corresponds to the customs of Greek
mathematics and is a consequence of the fact that the Greeks lacked
the technique of mathematical induction.

The essential weakness of IX.14 compared with the theorem of unique
prime factorization is the fact that IX.14 yields unique prime
factorization only for products of distinct primes, i.e., for
squarefree numbers; in the end, for a product
$n = p_1^{a_1} p_2^{a_2} p_3^{a_3}$, IX.14 only says that $p_1$,
$p_2$ and $p_3$ are the only prime factors of
$n$\footnote{It should be added, however, that elsewhere Euclid
proves theorems that could be used to sharpen IX.14: thus it is
proved in Book IX, among other things (IX.13), that the only divisors
of a prime power $p^n$ are the numbers $1$, $p$, $p^2$, \ldots,
$p^{n-1}$ (for Euclid, no number divides itself). Finally, Euclid
also proved (IX.36) that the only divisors of the number $2^{p-1}q$
for prime $q = 2^p-1$ are given by $1$, $2$, $4$, \ldots, $2^{p-1}$,
$q$, $2q$, $4q$, \ldots, $2^{p-2}q$.}.

Despite this shortcoming, Euclid's IX.14, just like the theorem on
unique prime factorization today, is proved using ``Euclid's
Lemma'' (Prop.\ VII.30), according to which $p \mid ab$ for
prime\footnote{Euclid's notion of a prime number coincides with what
is called irreducible in modern algebra. The property in Euclid's
Lemma is then used to define primes.} $p$ always implies $p \mid a$
or $p \mid b$. What kept Euclid from a complete formulation and proof
of the fundamental theorem of arithmetic was, on the one hand,
certainly the inadequate language he had to use, and on the other
hand the fact that IX.14 was anything but fundamental for Euclid,
being not much more than a small lemma that he mentions in passing.
Euclid does not use his Proposition IX.14 anywhere afterwards, and
this observation alone forbids us to speak of a fundamental theorem
in this context.

The foundation of Euclidean number theory is the existence of the
Euclidean algorithm, and the fundamental theorem from which the
central results are derived is Euclid's Proposition VII.19, to which
we now turn.

\section{Euclid's Proposition VII.19}

In order to understand Euclid's Proposition VII.19 we need his
Definition VII.20:
\begin{quote}
  \emph{Numbers are proportional when the first is the same multiple,
    or the same part, or the same parts, of the second that the third
    is of the fourth.}
\end{quote}
That the first number is the same multiple of the second as the third
is of the fourth means that there is a natural number $m$ with
$a = mb$ and $c = md$. Similarly, ``the same part'' means that there
is a divisor $n$ of $b$ and $d$ with $a = \frac1n b$ and
$c = \frac1n d$; with $b = nx$ and $d = ny$ this gives $a = x$,
$b = nx$, $c = y$, $d = ny$. Finally, ``the same parts'' means
$a = \frac mn b$ and $c = \frac mn d$ for natural numbers $m, n$ with
$n \mid b$ and $n \mid d$.\footnote{We would regard the first two
cases as the special cases $n = 1$ and $m = 1$, respectively, of the
last definition; Euclid's concept of number, however, differs from
ours: for him, numbers were proper multiples of the unit $1$, and the
unit was not a number.} Translated, this means:

\medskip\noindent{\bf Definition.}
For natural numbers $a, b, c, d \in \N$ the ratios $a:b$ and $c:d$
are proportional (i.e., $a:b = c:d$) if and only if there exist
$m, n, x, y \in \N$ with
\begin{equation}\label{eEamx}
 a = mx, \quad b = nx, \quad c = my, \quad d = ny.
\end{equation}

Thus, for example, $6:2 = 9:3$ via $m = 3$, $n = 1$, $x = 2$ and
$y = 3$. Now we can state
\begin{quote}
  \emph{Proposition VII.19. If four numbers be proportional, the
    number produced from the first and fourth will be equal to the
    number produced from the second and third; and, if the number
    produced from the first and fourth be equal to that produced from
    the second and third, the four numbers will be proportional.}
\end{quote}
Thus Euclid claims that $a:b = c:d$ if and only if $ad = bc$. In the
subsequent Proposition VII.20 Euclid shows that the least numbers
$a, b$ having the ratio $c:d$ divide these numbers $c$ and $d$,
namely in such a way that $c = am$ and $d = bm$. Finally, VII.21 and
VII.22 show that $a$ and $b$ are the least numbers with
$a:b = c:d$ if and only if $a$ and $b$ are coprime.

We tend to identify Euclid's proportions of numbers with modern
fractions (for Euclid, proportions were not numbers, and he does not
explain how to add proportions); in this language, Propositions
VII.19\,ff.\ are concerned with showing that every fraction can be
written uniquely in lowest terms.

\subsection*{Of numbers and magnitudes}
Besides proportions of numbers, Euclid also treats proportions of
magnitudes in Book V; these include angles, lengths, areas and
volumes, i.e., objects that can vary continuously. All magnitudes of
the same kind in Euclid have the following properties: they can be
added, and they behave archimedeanly in the following sense: if $A$
and $B$ are magnitudes of the same kind (such magnitudes can be
compared with respect to their size), then there always exists an
$n \in \N$ such that $nB > A$.

For ratios $A:B$ and $C:D$ of magnitudes\footnote{Here, of course,
$A$ and $B$, as well as $C$ and $D$, must be magnitudes of the same
kind.}, a definition attributed to Eudoxus (Def.~V.5) says that
$A:B = C:D$ if and only if the following holds: for all natural
numbers $m, n$ we have $mA < nB$ ($mA = nB$; $mA > nB$) if and only
if $mC < nD$ ($mC = nD$; $mC > nD$).

With this one can show, for example, that similar polygons inscribed
in circles are to one another as the squares on the diameters of the
circles (Prop.\ XII.1), or that circles are to one another as the
squares on their diameters (Prop.\ XII.2).

In order to make the differences between the definitions of
proportionality for numbers and for magnitudes clear, we observe the
following:
\begin{itemize}
\item Euclid's definition of proportionality of numbers makes sense
      not only in $\N$, but in arbitrary multiplicatively written
      monoids.
\item The definition of proportionality of magnitudes makes sense
      for all additively written monoids that are ``archimedean'' in
      the above sense.
\end{itemize}

In principle, numbers could also be regarded as
magnitudes\footnote{After all, $\N$ is also an additive archimedean
monoid.}; Euclid, however, does not do this. The claim one sometimes
encounters that Euclid consciously or unconsciously applied
statements about proportions of magnitudes to numbers seems to me to
be largely unfounded, and it does not fit well with the fact that the
theory of proportions of numbers is considerably
older\footnote{In \cite{vdW}, van der Waerden attributes practically
all of Book VII to the Pythagoreans; as we will explain below, the
definition of proportionality of numbers almost suggests itself if,
like the Pythagoreans, one is concerned with harmonics and whole
numbers. The definition of proportionality of magnitudes, on the
other hand, is generally attributed to Eudoxus.} than that of
magnitudes.

\section{Consequences of Euclid's VII.19}
At first sight, Euclid's VII.19 is an inconspicuous result; in fact,
however, it is quite an extraordinary one. In what follows we will
briefly discuss some well-known consequences of VII.19. We begin with

\begin{prop}
Euclid's Lemma follows from VII.19: all irreducible elements are
prime.
\end{prop}

\begin{proof}
  Let $p$ be irreducible and $p \mid bc$, i.e., $pd = bc$. By VII.19
  there exist $m, n, x, y$ with $p = mx$, $b = nx$, $c = my$ and $d =
  ny$. Since $p$ is irreducible, $m$ or $x$ is a unit, and this
  implies $p \mid b$ in the first and $p \mid c$ in the second case.
\end{proof}

VII.19 also suffices to prove the following generalization of
Euclid's Lemma:

\begin{prop}\label{eP1}
Euclid's Proposition VII.19 for natural numbers is equivalent to the
following result: if $a \mid bc$ and $\gcd(a,b) = 1$, then
$a \mid c$.
\end{prop}

\begin{proof}
  Assume VII.19 holds. If $a \mid bc$, i.e., $ad = bc$ for some $d \in
  \N$, then there exist $m, n, x, y \in \N$ with $a = mx$, $b = nx$,
  $c = my$ and $d = ny$. Since $\gcd(a,b) = 1$ we must have $x = 1$,
  hence $a = m$ is a divisor of $c = my$.

  Now let $ad = bc$, and put $x = \gcd(a,b)$. Then there exist $m, n
  \in \N$ with $a = mx$ and $b = nx$. Hence $md = nc$ and $\gcd(m,n) =
  1$; this implies $m \mid c$, i.e., $c = my$ for some $y \in \N$. But
  then $d = \frac{nc}m = ny$, i.e., (\ref{eEamx}) holds.
\end{proof}

Nowadays it is of course no problem to prove the uniqueness of prime
factorization with the help of Prop.\ \ref{eP1} (or Euclid's Lemma).
Using an idea of Zermelo \cite{Zer}, however, one can derive it
directly from Euclid's VII.19:

\begin{prop}
The uniqueness of factorization into irreducible elements follows
from Euclid's VII.19.
\end{prop}

\begin{proof}
  If uniqueness were false, there would be a smallest natural number
  $N$ with two different factorizations. Let $p$ be an irreducible
  factor of the first factorization and write $N = ph$; similarly, let
  $q$ be an irreducible factor of the second factorization, and $N =
  qk$. From $ph = qk$, VII.19 gives the existence of $m, n, x, y \in
  \N$ with $p = mx$, $q = nx$, $k = my$ and $h = ny$.  Since $p$ is
  irreducible, we must have $m = 1$ or $x = 1$. But then $p \mid q$
  (and hence $p = q$) in the first and $p \mid k$ in the second
  case. In the second case we have $k < N$, so by the choice of $N$
  the number $k$ has a unique factorization into irreducible elements;
  since $p \mid k$, the factor $p$ occurs in it, and hence also in the
  factorization of $k$ obtained from the second factorization of $N$
  by omitting $q$. Thus in both cases $p$ also occurs in the second
  factorization, and cancelling $p$ yields a natural number $M = \frac
  Np < N$ with two different factorizations into irreducible elements:
  this contradiction to the minimality of $N$ proves the claim.
\end{proof}

Kalm\'ar \cite{Kalm} (cf.\ Sur\'anyi \cite{Sur1,Sur2}) calls
Euclid's VII.19 the four number theorem and attributes it to Euler.
The generalization to more numbers (starting from a factorization
$a_1 \cdots a_i = b_1 \cdots b_j$) was also already known to Euler;
it follows from the simple four number theorem by induction on the
number of factors. Kalm\'ar then observed that the uniqueness of
prime factorization follows from this generalization (see also Bell
\cite{Bel1,Bel2}, Rosenthall \cite{Ros} and Erd\H{o}s \& Sur\'anyi
\cite{ES}), and Sur\'anyi \cite[p. 44\,ff.]{Sur2} gave a geometric
proof of the four number theorem. It seems that nobody has so far
noticed that the four number theorem is of Euclidean
origin\footnote{Anyone who has studied the Elements more closely
could of course have said so immediately, but apparently there are
rather few mathematicians who are sufficiently familiar both with
the Elements and with modern algebra. Even Sur\'anyi, who in
\cite{Sur2} deals quite explicitly with Greek mathematics, overlooked
the connection between Euclid's VII.19 and the four number
theorem.}.

It is an easy exercise to prove, e.g.\ with the help of the theorem
of unique prime factorization, the following result, which is used
in many elementary investigations of diophantine equations:

\begin{prop}\label{ePsq}
If $a, b \in \N$ are coprime and $ab$ is a square, then $a$ and $b$
are also squares.
\end{prop}

This proposition, too, is a simple consequence of Euclid's VII.19,
as is easily verified (we will come back to the proof in
Prop.\ \ref{ePsqc}). Gauss was the first to point out in the
Disquisitiones \cite{DA} that this result (or its generalization to
$n$-th powers) follows from unique prime factorization in $\Z$. The
converse conclusion, that Prop.\ \ref{ePsq} can only be proved with
the help of the fundamental theorem of arithmetic, is often
encountered, but false.

Euler, of whom it is often said that he used unique prime
factorization unconsciously in his works or regarded it as obvious,
mentioned Prop.\ \ref{ePsq} in his first paper \cite{EuTQ} on
diophantine questions:
\begin{quote}
  \emph{Factum ex duobus pluribusque numeris inter se primis nec
    quadratum nec cubus nec ulla alia potestas esse potest, nisi
    singuli factores sint quadrata vel cubi vel eiusmodi aliae
    potestates.\footnote{A product of two or more pairwise coprime
      numbers cannot be a square, nor a cube, nor any other power,
      unless the individual factors are squares or cubes or other
      powers of the same kind.}}
\end{quote}
Instead of giving a proof, he then remarks:
\begin{quote}
  \emph{Demonstratio huius Lemmatis facilis est atque ab Euclide iam
    est tradita, ita ut superfluum foret eam hic
    exponere.\footnote{The proof of this lemma is easy and has already
      been handed down by Euclid, so that it would be superfluous to
      present it here.}}
\end{quote}
I do not know which of Euclid's results Euler is referring to here.
For squares, the claim follows relatively easily from Euclid's theory
of similar plane numbers (in modern terminology, similar plane
numbers are numbers whose product is a square; cf.\ IX.1 and IX.2),
and for higher powers, with somewhat more effort, from many of his
propositions in Book VIII.

Thus where we would refer to unique prime factorization, i.e., to an
introductory course in number theory or to Gauss's Disquisitiones
\cite{DA}, Euler settles the matter with a reference to Euclid. It
seems quite plausible to me that in other places as well, where we
see a gap that can be closed by unique prime factorization, Euler
would have dismissed any objections with a reference to Euclid.

In \cite{Eu152}, for example, he introduces, for the study of
amicable numbers, the sum $\sigma(n)$ of all divisors of $n \in \N$
and shows in Lemma 1 that $\sigma$ is multiplicative, i.e., that
$\sigma(mn) = \sigma(m)\sigma(n)$ for coprime numbers $m, n$. Here is
his proof: ``For the product $mn$ has, first, the individual
divisors of both factors $m$ and $n$, and is moreover divisible by
the products of the divisors of $m$ and those of $n$. Therefore the
sum of all divisors of $mn$ is equal to the product of $\sigma(m)$
and $\sigma(n)$.'' Thus Euler uses the following lemma without proof:

\begin{lem}\label{eLEu4}
Let $m$ and $n$ be coprime. Then every divisor $a$ of $mn$ can be
written uniquely in the form $a = rs$ with $r \mid m$ and $s \mid n$.
\end{lem}

Today we would derive this lemma without hesitation from the
fundamental theorem of arithmetic; Euler, on the other hand, would
presumably have referred to Euclid: the existence of $r$ and $s$
follows like this: $a \mid mn$ means $ab = mn$ for some $b \in \N$;
Euclid's VII.19 yields the existence of $r, s, t, u \in \N$ with
$a = rs$, $b = tu$, $m = rt$ and $n = su$, and the claim follows.
Uniqueness is clear as well: if $a = rs = r's'$ with $r, r' \mid m$
and $s, s' \mid n$, then $r \mid r's'$ and
$\gcd(r,s') \mid \gcd(m,n) = 1$, hence $r \mid r'$ by
Prop.\ \ref{eP1}. Similarly $r' \mid r$, so $r = r'$ and $s = s'$.

\section{The gap}
As we have seen, the validity of Euclid's VII.19 is responsible for
the uniqueness of factorization into irreducible elements in $\N$. A
simple structure in which VII.19 no longer holds is the monoid
$M = \{1, 5, 9, 13, \ldots\}$ of all natural numbers
$\equiv 1 \bmod 4$ (in his lectures \cite{HZ}, Hilbert used the
monoid of all natural numbers $\equiv 1 \bmod 5$ in order to motivate
the introduction of ideals). In $M$ we clearly have
$21 \cdot 21 = 9 \cdot 49$ (this equation also shows immediately
that the analogue of Prop.\ \ref{ePsq} does not hold in $M$), but
$49:21 \ne 21:9$; for (\ref{eEamx}) would imply the existence of
$m, n, x, y \in M$ with $49 = mx$, $21 = nx = my$ and $9 = ny$. Since
$9$ is irreducible in $M$, we must have $n = 1$ or $y = 1$. In the
first case $y = 9$, hence $9 \mid my = 21$; in the second case
$n = 9$, and again $9 \mid nx = 21$; both are impossible.

On the other hand, one easily shows that $a:b = af:bf$ always holds
for $a, b, f \in M$. Hence
$$ 49 : 21 = 49\cdot 9 : 21 \cdot 9
           = 21 \cdot 21 : 21 \cdot 9 = 21 : 9, $$
where we have applied the observation $a:b = af:bf$ twice. But as we
have already seen, the equation $49 : 21 = 21 : 9$ is false. Where
did we make the mistake? The equations
$49:21 = 21 \cdot 21 : 21 \cdot 9$ and
$21 \cdot 21 : 21 \cdot 9 = 21 : 9$ are, as one can also easily check
directly, undoubtedly correct. The error in the argument above lies
in the assumption that equality of proportions\footnote{One might be
tempted to attribute this to the fact that the transitivity of
equality is an axiom in Euclid (``Things which are equal to the same
thing are also equal to one another''); however, what we call
equality of proportions is not an equality in the Euclidean sense:
Euclid does not define the equality of two proportions $a:b$ and
$c:d$, but says that four numbers are proportional if there exist
numbers as in (\ref{eEamx}). For this reason, many commentators of
Euclid do not write $a:b = c:d$, but more precisely $a:b :: c:d$.}
is transitive: from $a:b = c:d$ and $c:d = e:f$ it does not follow in
general that $a:b = e:f$.

Taisbak discovered in \cite{Tai} that Euclid left exactly this gap in
the proof of his fundamental theorem VII.19: he used the
transitivity of equality, but did not prove it. How to close this gap
can be found in Taisbak \cite{Tai} as well as in Pengelley \&
Richman \cite{PR}. The simple argument is based on showing that some
of the numbers in the four number theorem can be chosen more or less
``canonically'' as greatest common divisors:

\begin{lem}\label{eLtr}
Let $a:b = c:d$. Then there exist $p, q, r, s \in \N$ with
$a = pr$, $b = qr$, $c = ps$ and $d = qs$, where we may choose
$r = \gcd(a,b)$ and $s = \gcd(c,d)$.
\end{lem}

\begin{proof}
  From $a:b = c:d$ we get the existence of numbers $m, n, x, y \in \N$
  with (\ref{eEamx}). Thus $x \mid \gcd(a,b)$ and $y \mid \gcd(c,d)$,
  hence there exist $i, j \in \N$ with $\gcd(a,b) = ix$ and $\gcd(c,d)
  = jy$. From $ix \mid a = mx$ and $ix \mid b = nx$ we obtain $i \mid
  m$ and $i \mid n$, i.e., $m = pi$ and $n = qi$ with $p, q \in
  \N$. In particular, $a = pix = p \cdot \gcd(a,b)$ and $b = q \cdot
  \gcd(a,b)$.

  We now claim that $i \mid j$. By symmetry (the situation is
  symmetric under exchanging $(a,b,x)$ and $(c,d,y)$) we then also
  have $j \mid i$, hence $i = j$. But then we get, as desired,
  $c = my = piy = pjy = p \cdot \gcd(c,d)$ and
  $d = ny = qiy = qjy = q \cdot \gcd(c,d)$.

  It remains to show $i \mid j$. From $m = pi$ and $n = qi$ we get $iy
  \mid my = c$ and $iy \mid ny = d$. Hence $iy \mid \gcd(c,d) = jy$,
  and therefore $i \mid j$ as claimed.
\end{proof}

With this we can prove the transitivity of equality as follows: let
$a:b = c:d$ and $c:d = e:f$. By Lemma \ref{eLtr} there exist
$p, q \in \N$ with $a = p \cdot \gcd(a,b)$, $b = q \cdot \gcd(a,b)$,
$c = p \cdot \gcd(c,d)$ and $d = q \cdot \gcd(c,d)$. Applying the
lemma a second time we find $p', q' \in \N$ with
$c = p' \cdot \gcd(c,d)$, $d = q' \cdot \gcd(c,d)$,
$e = p' \cdot \gcd(e,f)$ and $f = q' \cdot \gcd(e,f)$. Cancelling
gives $p = p'$ and $q = q'$, and then certainly $a:b = e:f$.

Note that Lemma \ref{eLtr} does not hold in the monoid $M$ above,
since in $M$ two numbers need not have a gcd. For example, the common
divisors of $3^2 \cdot 7 \cdot 11$ and $3^2 \cdot 7 \cdot 19$ in $M$
are $1$, $3^2$ and $3 \cdot 7$, but there is no greatest common
divisor with respect to divisibility, i.e., no common divisor that is
divisible by all common divisors.

\section{Euclid or Pythagoras?}

Numerous historians of mathematics attribute Euclid's approach to
number theory to the Pythagoreans. They are assumed to have arrived
at their proportions by way of music: they had observed that the
lengths of vibrating strings stand in simple ratios to one another
when the corresponding tones sound harmonious. Thus the ratios $1:2$
and $2:3$ correspond to the octave and the fifth, respectively.

One of the first observations must then have been that certain
proportions can be simplified: after all, $2:4 = 1:2$. Ultimately,
this line of inquiry inevitably leads to the result that proportions
can be reduced uniquely: if $m$ and $n$ are minimal with
$a:b = m:n$, then $m \mid a$ and $n \mid b$; more precisely, there
exists an $x \in \N$ with $a = mx$ and $b = nx$. But this is the
content of Euclid's Prop.\ VII.20. Afterwards, VII.22 shows that the
minimal $m$ and $n$ are coprime.

Once VII.20 has been discovered empirically, and if $a:b = c:d$, it
follows on the one hand that $a = mx$ and $b = nx$, and on the other
hand that $c = my$ and $d = ny$. But this is exactly Euclid's
definition of proportionality of numbers, which thus arises quite
naturally from the study of the most elementary properties of
proportions.

According to van der Waerden \cite{vdW}, the Pythagoreans made this
four number theorem the basis of their number theory, which Euclid
later took over, together with its structure, essentially unchanged.
Like other commentators before and after him, van der Waerden
regarded Euclid's definition of proportionality of numbers as a
blunder\footnote{D. Fowler, for example, writes in \cite{Fow}:
\begin{quote}
  \emph{Commentators rarely point out the unsatisfactory nature of
    this definition: it is a vivid, though incomplete, description of
    four numbers in proportion, and not the mathematical criterion
    needed for the foundation of a theory.}
\end{quote}}
rather than as a stroke of genius. The above considerations and,
e.g., the proof of Lemma \ref{eLEu4} show, however, that the
definition is on the one hand quite natural, and that on the other
hand the four number theorem is perfectly suited for proving
fundamental results of diophantine analysis.

Besides the arguments given by van der Waerden for the claim that
Euclid's Book VII goes back to the Pythagoreans, one may also regard
the missing proof of the transitivity of equality for proportions of
numbers as evidence for it: if one bases the theory on the definition
of equality of proportions given by Eudoxus, the proof of the
transitivity of equality becomes almost trivial. Nevertheless,
Euclid carried it out in V.11. This may be due to the fact that the
theory of proportions of magnitudes goes back to Eudoxus, and that
in his time the insight that such steps have to be proved was
already widespread among mathematicians. In my opinion it is quite
conceivable that the Pythagoreans simply overlooked the fact that the
transitivity of equality of proportions of numbers is not trivial at
all. That Euclid did not see this either, although he of course knew
V.11, is more surprising, but could be explained by the fact that he
took over the number theory of the Pythagoreans into his Book VII
without major changes.

\section{The four number theorem}

In order to clarify the significance of Euclid's four number theorem,
we call integral domains in which this theorem holds Riesz
domains\footnote{Zafrullah \cite{Z} called these rings pre-Schreier;
below, however, we will see that the divisibility groups of Riesz
domains are ``Riesz groups'' in the sense of group theory.}. Thus in
a Riesz domain $R$, $ad = bc$ with $a, b, c, d \in R$ implies the
existence of $m, n, x, y \in R$ with $a = mx$, $b = nx$, $c = my$ and
$d = ny$.

An element $a \in R$ is called primal\footnote{This notion is due to
Cohn \cite{Cohn}.} if $a \mid bc$ for arbitrary $b, c \in R$ implies
the existence of $r, s \in R$ with $a = rs$, $r \mid b$ and
$s \mid c$. Note that $0$ and the units are always primal.

Clearly an integral domain is a Riesz domain if and only if every
element of $R$ is primal. The following proposition shows that in
Riesz domains all irreducible elements are prime:

\begin{prop}\label{eLIP}
Irreducible elements are primal if and only if they are prime.
\end{prop}

\begin{proof}
  Let $a$ be prime and $a \mid cd$; then without loss of generality $a
  \mid c$. With $r = a$ and $s = 1$ we then have $a = rs$, $r \mid c$
  and $s \mid d$, i.e., $a$ is primal.

  Now let $a$ be irreducible and primal. If $a \mid cd$, there exist
  $r, s$ with $a = rs$, $r \mid c$ and $s \mid d$. Since $a$ is
  irreducible, $r$ or $s$ -- say $r$ -- must be a unit. But then $s
  \mid d$ implies $a \mid d$, and so $a$ is prime.
\end{proof}

In Prop.\ VII.19, Euclid essentially showed that $\Z$ is a Riesz
domain; his proof that irreducible elements in $\N$ are prime is
essentially one half of the proof above.

A closer look at Euclid's proof of VII.19 (or at its completion by
Taisbak) shows that the four number theorem holds in all GCD domains
(these are integral domains in which any two elements have a greatest
common divisor), i.e., that every GCD domain is a Riesz domain.

At this point it should be pointed out that primal elements on the
one hand and irreducible or prime elements on the other hand are of
quite different nature: products of several irreducible (or prime)
elements are no longer irreducible (or prime), apart from the trivial
case of a single factor. For primal elements, on the other hand, we
have:

\begin{lem}
Products of primal elements are primal.
\end{lem}

\begin{proof}
  Let $a, b \in R$ be primal; if $a = 0$ or $b = 0$, then $ab = 0$ is
  primal, so we may assume $a, b \ne 0$. Let $ab \mid c_1c_2$, i.e.,
  $c_1c_2 = abd$. Since $a$ is primal, there is a factorization
  $a = a_1a_2$ with $a_i \mid c_i$. Writing $c_i = a_id_i$, we get
  $abd = c_1c_2 = a_1a_2d_1d_2 = ad_1d_2$, hence $bd = d_1d_2$. Since
  $b$ is primal, we have $b = b_1b_2$ with $b_i \mid d_i$.

  Thus $ab = a_1b_1a_2b_2$ with $a_ib_i \mid a_id_i = c_i$, and
  consequently $ab$ is primal. By induction it follows that arbitrary
  products of primal elements are primal.
\end{proof}

We now prove the following generalization of Prop.\ \ref{ePsq}:

\begin{prop}\label{ePsqc}
Let $R$ be a Riesz domain. If $ab = c^2$ for $a, b, c \in R$ and if
$a$ and $b$ are coprime, then there is a unit $e \in R^\times$ such
that $ae$ and $be$ are squares in $R$.
\end{prop}

\medskip\noindent
{\bf Remark.}
Although the notion of a greatest common divisor of two arbitrary
ring elements only makes sense in GCD domains, one may speak of
``coprime'' elements in arbitrary integral domains: by definition,
these are elements whose common divisors are all units. For example,
$2$ and $3$ are coprime in $R = \Z[\sqrt{-5}\,]$, whereas $6$ and
$2+2\sqrt{-5}$ do not have a gcd in $R$.

\medskip

\begin{proof}[Proof of Prop.\ \ref{ePsqc}]
  By the four number theorem there exist $m, n, x, y \in R$ with $a =
  mx$, $b = ny$ and $c = my = nx$. Applying the four number theorem to
  the equation $my = nx$ we get the existence of $r, s, t, u \in R$
  with $m = rs$, $y = tu$, $n = rt$ and $x = su$.  Since $a$ and $b$
  are coprime, so are $m$ and $n$, as well as $x$ and $y$;
  consequently $r$ and $u$ are units in $R$. With $e = 1/ru$ we then
  have $ae = mxe = rs^2ue = s^2$ and $be = nye = rt^2ue = t^2$.
\end{proof}

Even though some facts used by Euler without proof can be justified
with Euclid's four number theorem, Riesz domains are only ``almost''
factorial. If we consider two factorizations
$a_1 \cdots a_m = b_1 \cdots b_n$, then in Riesz domains it only
follows that both factorizations have a common refinement: there
exist elements $c_{ij} \in R$ ($1 \le i \le m$, $1 \le j \le n$)
with $a_i = c_{i1} \cdots c_{in}$ and $b_j = c_{1j} \cdots c_{mj}$.
For $m = n = 2$ this is just the four number theorem; the general
case (which follows from the four number theorem by induction) is
called the ``Schreier refinement property''. For example, the ring
of all algebraic integers has this refinement property, but it is
not factorial since it has no irreducible elements.

Integral domains in which every nonzero element can be written as a
product of a unit and irreducible elements are called atomic. Thus
we have the following

\begin{thm}
An integral domain $R$ is factorial if and only if $R$ is an atomic
Riesz domain.
\end{thm}

\begin{proof}
  In Riesz domains, irreducible elements are prime by
  Prop.\ \ref{eLIP}; as is well known, this implies that
  factorizations into irreducible elements are unique. Conversely, if
  $R$ is factorial, then $R$ is a GCD domain and hence a Riesz domain,
  and secondly every element $\ne 0$ has a factorization into
  irreducible (even prime) elements.
\end{proof}

Because of Prop.\ \ref{eLIP}, this statement is somewhat weaker than
the well-known theorem according to which integral domains are
factorial if and only if they are atomic and every irreducible
element is prime.

Euclid's construction of number theory can be schematized as
follows:
$$ \text{Euclid} \lra \text{GCD} \lra \text{Riesz} \lra \text{AP}, $$
where AP stands for the validity of Euclid's Lemma, i.e., for the
property that irreducible elements (atoms) are prime. More generally,
one can show the following inclusions:
\begin{figure}[h!]
$$ \begin{array}{cccccccccccc}
\text{Euclid} & \subset &  \text{PID}    & \subset & \text{Bezout}
        & & & & & & \\
& &     \cap      &         &   \cap     &         & & & & & \\
& & \text{UFD} & \subset & \text{GCD} & \subset & \text{Riesz}
              & \subset & \text{AP} \\
& &     \cap      & &       &    & & & & & \\
& & \text{atomic}  & &       &    & & & & &
\end{array}$$
\end{figure}

Here Euclid denotes Euclidean domains, PID principal ideal domains,
UFD factorial domains (i.e., those with unique factorization into
prime elements), and AP domains in which all atoms (i.e., all
irreducible elements) are prime. Bezout domains are those in which
any two elements $a, b$ have a gcd that can be written as an
$R$-linear combination of $a$ and $b$; equivalently, the sum of two
principal ideals is again a principal ideal.

While Hilbert put Euclidean geometry on a modern axiomatic
foundation, a corresponding investigation of Euclidean number theory
has, to my knowledge, not yet been carried out. One may ask, for
example, in which rings the following holds: if $a$ and $b$ are
coprime (i.e., every common divisor is a unit) and $a \mid bc$, then
$a \mid c$. The first part of the proof of Prop.\ \ref{eP1} shows
that this holds in Riesz domains; the second part of the proof,
however, uses the existence of the gcd, so that we cannot conclude
that such rings are also Riesz domains. On the other hand, rings with
this property are AP. In general one may ask in which integral
domains the propositions of Euclid's number theory hold, i.e., which
axioms are needed to prove them.

\section{Riesz groups}

In this section we will discuss some generalizations of the four
number theorem. To this end we will study, somewhat more generally,
Riesz groups.

Let $R$ be an integral domain with quotient field $K$. Recall that
$a \mid c$ for $a, c \in K^\times$ means that there is a $b \in R$
with $ab = c$. Units are irrelevant for the study of divisibility
properties; it therefore makes sense to consider, instead of
$K^\times$, the factor group $G(R) = K^\times/R^\times$, which we can
also regard as the group of nonzero fractional principal ideals; for
the time being we denote the elements of $G(R)$ by $aR^\times = (a)$.
The group $G(R)$ is partially ordered by $(a) \le (b)$ if and only if
$a \mid b$, i.e., if $\frac ba \in R$. Moreover, $(a) \le (b)$ is
equivalent to $(ac) \le (bc)$ for arbitrary $(c) \in G(R)$. The set
$G(R)^+ = \{(g) \in G(R): (g) \ge (1)\}$ of all ``positive'' elements
thus consists of all nonzero integral principal ideals and is clearly
a monoid, which arises from $R^* = R \setminus \{0\}$ by passing to
ideals, i.e., by identifying associated elements.

We now want to translate Euclid's four number theorem for $R$ into
this language. If $ad = bc$, then $a \mid bc$, hence $(a) \le (b)(c)$
for the corresponding principal ideals in $G(R)$. By the four number
theorem there exist $r, s, t, u \in R$ with $a = rs$, $b = ru$,
$c = ts$ and $d = tu$. Thus $r \mid b$ and $s \mid c$, and hence
$(r) \le (b)$ as well as $(s) \le (c)$. Therefore we have

\begin{prop}
Let $R$ be an integral domain. Then the four number theorem holds in
$R$ if and only if $G(R)$ has the following ``Fuchs interpolation
property'' $F_2$\footnote{From now on we denote the elements of
$G(R)$ simply by $a$ instead of $(a)$.}: if $a, b, c \in G(R)^+$
and $a \le bc$, then there is a factorization $a = rs$ with
$r, s \in G(R)^+$ such that $r \le b$ and $s \le c$.
\end{prop}

We will now define a series of general interpolation and refinement
properties of partially ordered abelian groups $G$:
\begin{itemize}
\item The ``Fuchs interpolation property'' $F_n$: if
      $a$, $b_1$, \ldots, $b_n \in G^+$ and $a \le b_1 \cdots b_n$,
      then there exist $a_i \in G^+$ with $a = a_1 \cdots a_n$ and
      $a_i \le b_i$ for $1 \le i \le n$.
\item The ``Riesz interpolation property'' $R_{m,n}$: if
      $a_1, \ldots, a_m$, $b_1, \ldots, b_n \in G$ are given with
      $a_i \le b_j$ for all $1 \le i \le m$ and $1 \le j \le n$, then
      there is a $c \in G$ with $a_i \le c \le b_j$ for all $i$ and
      $j$.
\item The ``Schreier refinement property'' $S_{m,n}$: if
      $a_1, \ldots, a_m$, $b_1, \ldots, b_n \in G^+$ are given with
      $a_1 \cdots a_m = b_1 \cdots b_n$, then there exist elements
      $c_{ij} \in G^+$ ($1 \le i \le m$, $1 \le j \le n$) with
      $a_i = c_{i1} \cdots c_{in}$ and $b_j = c_{1j} \cdots c_{mj}$.
\item The ``Euclidean refinement property'' $E_{m,n}$: if
      $a, b_i, c_j \in G^+$ and $a \le b_ic_j$ for $i = 1$, \ldots,
      $m$ and $j = 1$, \ldots, $n$, then there is a factorization
      $a = rs$ with $r, s \in G^+$, $r \le b_i$ and $s \le c_j$ for
      all $1 \le i \le m$ and $1 \le j \le n$.
\end{itemize}

The properties $R_{2,0}$ and $R_{0,2}$ can be interpreted as
follows: for $a_1, a_2 \in G$ there is a $c \in G$ with $a_1 \le c$
and $a_2 \le c$ (this is $R_{2,0}$), and a $b \in G$ with
$b \le a_1$ and $b \le a_2$ (this is $R_{0,2}$); partially ordered
groups with this property are called directed\footnote{Since
$a \le b$ is equivalent to $b^{-1} \le a^{-1}$, a partially ordered
group is already directed if, e.g., $R_{2,0}$ is satisfied.}. The
divisibility groups $G(R)$ of integral domains are automatically
directed: for $a_1, a_2 \in K^\times$ there is, as one sees
immediately, always a $c \in K^\times$ with $a_1 \mid c$,
$a_2 \mid c$ (choose $c$ as the product of the (arbitrarily chosen)
numerators of $a_1$ and $a_2$). Similarly one always finds a
$b \in K^\times$ with $b \mid a_1$ and $b \mid a_2$.

We now claim

\begin{thm}\label{eSVerf}
For a partially ordered directed abelian group $G$ the following
properties are equivalent:
\begin{enumerate}
\item\label{eR1} $R_{2,2}$ holds in $G$;
\item\label{eR2} $R_{m,n}$ holds in $G$ for all $m, n \ge 0$;
\item\label{eR3} $F_2$ holds in $G$;
\item\label{eR4} $F_n$ holds in $G$ for all $n \ge 1$;
\item\label{eR5} $S_{2,2}$ holds in $G$;
\item\label{eR6} $S_{m,n}$ holds in $G$ for all $m, n \ge 1$;
\item\label{eR7} $E_{1,1}$ holds in $G$;
\item\label{eR8} $E_{m,n}$ holds in $G$ for all $m, n \ge 1$;
\item\label{eR9} $[1,a] \cdot [1,b] = [1,ab]$ holds in $G$ for all
       $a, b \in G^+$ (here $[a,b] = \{r \in G: a \le r \le b\}$).
\end{enumerate}
\end{thm}

The assumption that $G$ is directed is only needed in order to pass
from $R_{2,2}$ to $R_{m,n}$. The proofs are all elementary
(essentially induction) and need not be given here.

Following L.~Fuchs \cite{Fu}, partially ordered abelian groups having
the properties listed in Theorem \ref{eSVerf} are called Riesz
groups\footnote{Strictly speaking, commutativity is not required for
Riesz groups; in this case, not all properties of
Theorem \ref{eSVerf} are equivalent anymore.}. They were first studied
by F.~Riesz \cite{Rie} in 1940, later by Birkhoff \cite{Birkh} and
H.~Bauer \cite{Bauer}. In these papers one also finds parts of the
proof of Theorem \ref{eSVerf} and further properties equivalent to
\ref{eR1} -- \ref{eR9}. The Schreier refinement property $S_{m,n}$ is
the generalization of the four number theorem that Kalm\'ar used to
prove the fundamental theorem of arithmetic. In the second part we
will finally need the Euclidean refinement property.

\section*{Concluding remarks}

Anyone who studies Euclid's number theory more closely will quickly
come across the central role of his Proposition VII.19, the four
number theorem. But anyone who is merely interested in the history of
number theory before Gauss must also know Euclid's number-theoretic
books, and in particular the central result VII.19: the frequently
voiced criticism that Euler in some places implicitly used a result
of Gauss, namely unique prime factorization, is, first, somewhat
anachronistic in this form anyway and, second, probably untenable in
view of Euclid's VII.19\footnote{Other authors (even after Gauss)
also used the four number theorem without referring to Euler or
Euclid: see Barlow \cite[p. 401]{Bar11}, Kramer \cite{Kramer}, or
P\'epin \cite[p. 24]{PepLeg}, \cite[p. 413]{Pep752}.}. These remarks
are not meant to absolve Euler of all ``guilt'': a considerable
number of his proofs are not rigorous by today's standards, and not
all gaps can be filled with Euclid. The sweeping accusation,
however, that Euler often tacitly used unique prime factorization
should perhaps be reconsidered.

\bigskip

This article was published in Math. Semesterberichte \textbf{55}
(2008), 181--195. The English translation was done with the help of
Claude (anthropic).

\end{otherlanguage*}

\end{document}